\magnification=1200
\catcode`\@=11
\font\tenmsa=msam10
\font\sevenmsa=msam7
\font\fivemsa=msam5
\font\tenmsb=msbm10
\font\sevenmsb=msbm7
\font\fivemsb=msbm5
\newfam\msafam
\newfam\msbfam
\textfont\msafam=\tenmsa  \scriptfont\msafam=\sevenmsa
  \scriptscriptfont\msafam=\fivemsa
\textfont\msbfam=\tenmsb  \scriptfont\msbfam=\sevenmsb
  \scriptscriptfont\msbfam=\fivemsb
\def\hexnumber@#1{\ifcase#1 0\or1\or2\or3\or4\or5\or6\or7\or8\or9\or
        A\or B\or C\or D\or E\or F\fi }
\edef\msa@{\hexnumber@\msafam}
\edef\msb@{\hexnumber@\msbfam}
\mathchardef\boxtimes="0\msa@02
\mathchardef\square="0\msa@03
\mathchardef\subsetneq="3\msb@28
\def\Bbb{\ifmmode\let\next\Bbb@\else
        \def\next{\errmessage{Use \string\Bbb\space only in math
mode}}\fi\next}
\def\Bbb@#1{{\Bbb@@{#1}}}
\def\Bbb@@#1{\fam\msbfam#1}
\catcode`\@=12
\font\titlefont=cmbx12
\font\namesfont=cmr12
\catcode`\@=11
\font\tenmib=cmmib10
\font\sevenmib=cmmib7
\font\fivemib=cmmib5
\newfam\mibfam
\textfont\mibfam=\tenmib  \scriptfont\mibfam=\sevenmib
  \scriptscriptfont\mibfam=\fivemib
\def\hexnumber@#1{\ifcase#1 0\or1\or2\or3\or4\or5\or6\or7\or8\or9\or
        A\or B\or C\or D\or E\or F\fi }
\edef\mib@{\hexnumber@\mibfam}
\def\bfit{\ifmmode\let\next\bfit@\else
        \def\next{\errmessage{Use \string\bfit\space only in math
mode}}\fi\next}
\def\bfit@#1{\fam\mibfam#1}
\font\teneufm=eufm10
\font\seveneufm=eufm7
\font\fiveeufm=eufm5
\newfam\eufmfam
\textfont\eufmfam=\teneufm  \scriptfont\eufmfam=\seveneufm
  \scriptscriptfont\eufmfam=\fiveeufm
\def\hexnumber@#1{\ifcase#1 0\or1\or2\or3\or4\or5\or6\or7\or8\or9\or
        A\or B\or C\or D\or E\or F\fi }
\edef\eufm@{\hexnumber@\eufmfam}
\def\frak{\ifmmode\let\next\frak@\else
        \def\next{\errmessage{Use \string\frak\space only in math
mode}}\fi\next}
\def\frak@#1{\fam\eufmfam#1}
\catcode`\@=12
\font\eightrm=cmr8  \font\sixrm=cmr6
\font\eighti=cmmi8  \font\sixi=cmmi6
\font\eightsy=cmsy8 \font\sixsy=cmsy6
\font\eightbf=cmbx8 \font\sixbf=cmbx6
\font\eighttt=cmtt8
\font\eightit=cmti8
\font\eightsl=cmsl8

\catcode`@=11
\newskip\ttglue
\def\eightpoint{\def\rm{\fam0\eightrm}
  \textfont0=\eightrm \scriptfont0=\sixrm \scriptscriptfont0=\fiverm
  \textfont1=\eighti  \scriptfont1=\sixi  \scriptscriptfont1=\fivei
  \textfont2=\eightsy \scriptfont2=\sixsy \scriptscriptfont2=\fivesy
  \textfont3=\tenex   \scriptfont3=\tenex \scriptscriptfont3=\tenex
  \textfont\itfam=\eightit  \def\it{\fam\itfam\eightit}%
  \textfont\slfam=\eightsl  \def\sl{\fam\slfam\eightsl}%
  \textfont\ttfam=\eighttt  \def\tt{\fam\ttfam\eighttt}%
  \textfont\bffam=\eightbf  \scriptfont\bffam=\sixbf
  \scriptscriptfont\bffam=\fivebf  \def\bf{\fam\bffam\eightbf}%
  \tt \ttglue=.5em plus.25em minus.15em
  \normalbaselineskip=9pt
  \setbox\strutbox=\hbox{\vrule height7pt depth2pt width0pt}%
  \let\sc=\sixrm \let\big=\eightbig \normalbaselines\rm}
\def\eightbig#1{{\hbox{$\textfont0=\ninerm\textfont2=\ninesy
  \left#1\vbox to6.5pt{}\right.\n@space$}}}
\catcode`@=12
\def\whiteA{\vskip 1.5pc plus.5pc minus.1pc}
\def\whiteB{\vskip 1.12pc plus.5pc minus.1pc}
\def\whiteC{\vskip.5pc plus.2pc minus.1pc}
\def\whiteD{\vskip.5pc plus.5pc minus.1pc}
\def\mC{{\Bbb C}}
\def\mG{{\Bbb G}}
\def\mQ{{\Bbb Q}}
\def\mR{{\Bbb R}}
\def\mS{{\Bbb S}}
\def\mZ{{\Bbb Z}}

\def\cB{{\cal B}}
\def\cD{{\cal D}}
\def\gc{{\frak c}}
\def\gg{{\frak g}}
\def\gh{{\frak h}}
\def\hg{{\frak hg}}
\def\gsl{{\frak sl}}
\def\uu{{\frak u}}
\def\mapright#1{\ \smash{\mathop{\longrightarrow}\limits^{#1}}\ }
\def\mapleft#1{\ \smash{\mathop{\longleftarrow}\limits^{#1}}\ }
\def\twoheadrightarrow{\rightarrow\kern -8pt\rightarrow}
\def\isomarrow{\mapright\sim}
\def\leftisomarrow{\mapleft\sim}
\def\Aut{{\rm Aut}}
\def\Cent{{\rm Cent}}
\def\End{{\rm End}}
\def\Gal{{\rm Gal}}
\def\GL{{\rm GL}}
\def\Hg{{\rm Hg}}
\def\Hom{{\rm Hom}}
\def\HS{{\rm HS}}
\def\Ker{{\rm Ker}}

\def\MT{{\rm MT}}
\def\Nm{{\rm Nm}}
\def\Res{{\rm Res}}
\def\SL{{\rm SL}}
\def\Sp{{\rm Sp}}
\def\SU{{\rm SU}}

\def\UU{{\rm U}}

\def\Ad{{\rm Ad}}

\def\det{{\rm det}}
\def\id{{\rm id}}
\def\opp{{\rm opp}}
\def\pr{{\rm pr}}
\def\trace{{\rm tr}}
\def\Romno#1{\uppercase\expandafter{\romannumeral #1}}

\newcount\chno
\def\chapno{\number\chno}
\newcount\secno
\def\sectno{\number\secno}
\let\tempcolon=\colon
\def\colon{\tempcolon\,}

\let\tempcirc=\circ
\def\circ{\hskip 0.1em%
\mathord{\raise0.25ex\hbox{$\scriptscriptstyle\tempcirc$}}\hskip 0.1em}

\let\tempboxtimes=\boxtimes
\def\boxtimes{\hskip 0.1em%
\mathord{\raise0.15ex\hbox{$\scriptstyle\tempboxtimes$}}\hskip 0.1em}

\let\phi=\varphi
\let\epsilon=\varepsilon

\def\nmid{\kern-1pt\not{\kern -.6pt |}\kern 2pt}

\def\Proof{{\sl Proof.\/}}
\def\gdot{{\raise.25ex\hbox{${\scriptscriptstyle \bullet}$}}}

\def\Eo{{\End^0}}
\newcount\cntr
\cntr=1
\def\incr{\advance\cntr by 1}
\edef\refBorov{\number\cntr}\incr
\edef\refBourLie{\number\cntr}\incr
\edef\refDelWeilK3{\number\cntr}\incr
\edef\refDelShim{\number\cntr}\incr
\edef\refHazCM{\number\cntr}\incr
\edef\refHazHC{\number\cntr}\incr
\edef\refHazNons{\number\cntr}\incr
\edef\refImai{\number\cntr}\incr
\edef\refJacobs{\number\cntr}\incr
\edef\refLZTate{\number\cntr}\incr
\edef\refMZDuke{\number\cntr}\incr
\edef\refMZCrel{\number\cntr}\incr
\edef\refMumNote{\number\cntr}\incr
\edef\refMAV{\number\cntr}\incr
\edef\refMurtyHC{\number\cntr}\incr
\edef\refOortEnd{\number\cntr}\incr
\edef\refPink{\number\cntr}\incr
\edef\refRibet{\number\cntr}\incr
\edef\refSchoenHC{\number\cntr}\incr
\edef\refSerRGF{\number\cntr}\incr
\edef\refSerGAHT{\number\cntr}\incr
\edef\refShimFam{\number\cntr}\incr
\edef\refTank{\number\cntr}\incr
\edef\refTankPrime{\number\cntr}\incr
\edef\refvGeemen{\number\cntr}\incr
\edef\refWeil{\number\cntr}\incr
\edef\refZarWts{\number\cntr}\incr
\edef\refZarLNM{\number\cntr}\incr
\edef\refZarTaFin{\number\cntr}\incr
\edef\refYZBM{\number\cntr}\incr
\noindent
\centerline{{\titlefont Hodge classes on abelian varieties of low dimension}}
\bigskip\bigskip
\centerline{{\namesfont B.~J.~J. Moonen\footnote*{{\rm Research made possible
by a fellowship of the Royal Netherlands Academy of Arts and Sciences.}}\quad
and\quad Yu.~G. Zarhin\footnote{**}{{\rm Supported by the National Science
Foundation.}}}}
\vskip 6pc

\centerline{{\bf Introduction.}}
\whiteA

\noindent
In this paper we study Hodge classes on complex abelian varieties $X$. If
$\dim(X) \leq 3$ then every Hodge class on $X$ is a linear combination of
products of divisor classes. (This is true for any smooth projective complex
variety $X$.) The same property holds true for self-products of simple abelian
varieties of prime dimension, as shown by Tankeev [\refTankPrime]; see also
Ribet's paper [\refRibet]. In [\refMZDuke] the authors showed that if $X$ is
simple of dimension 4 then every Hodge class is a linear combination of
products of divisor classes and Weil classes---if there are any. (The notion of
a Weil class shall be briefly reviewed in (1.9); for an elementary discussion
see also [\refYZBM].)

The aim of this note is to extend this to arbitrary abelian varieties of
dimension $\leq 5$. In order to state our main results, let us describe some
special cases. We start with dimension 4.
\whiteC

(a)~The abelian variety $X$ is isogenous to a product $X_1 \times X_2$ where
$X_1$ is an elliptic curve with complex multiplication by an imaginary
quadratic field $k$ and where $X_2$ is a simple abelian threefold such that
there exists an embedding $k \hookrightarrow \Eo(X_2)$.
\whiteC

(b)~The abelian variety $X$ is simple of dimension 4 such that $\Eo(X)$ is a
field containing an imaginary quadratic field $k$ which acts on the tangent
space $T_{X,0}$ with multiplicities $(2,2)$. (See \S 1 for further
explanation.)
\whiteC

(c)~The abelian variety $X$ is simple of dimension 4 with $D = \Eo(X)$ a
definite quaternion algebra over $\mQ$. (Type \Romno 3 in the Albert
classification.) Note that for every $\alpha \in D \setminus \mQ$ the
subalgebra $\mQ(\alpha) \subset D$ is an imaginary quadratic field.
\whiteC

(d)~The abelian variety $X$ is simple of dimension 4 with $\Eo(X) = \mQ$.
\whiteB

\advance\secno by 1
\edef\refMainThm{\number\chno.\number\secno}
\noindent
{\bf (\chapno.\sectno) Theorem.}\enspace {\sl Let $X$ be a complex abelian
variety with $\dim(X) \leq 4$. Write $V = H_1(X(\mC),\mQ)$ and let $\phi \colon
V \times V \rightarrow \mQ$ be the Riemann form associated to a polarization of
$X$. Write $D = \Eo(X)$ and let $\Sp_D(V,\phi)$ denote the centralizer of $D$
inside the symplectic group $\Sp(V,\phi)$.

{\rm (\romannumeral1)}~Suppose we are in case {\rm (a)} or {\rm (b)}. Then the
Hodge ring $\cB^\gdot(X)$ is generated by the subalgebra $\cD^\gdot(X)$ of
divisor classes together with the space of Weil classes $W_k \subset \cB^2(X)$.
The Hodge group $\Hg(X)$ is strictly contained in $\Sp_D(V,\phi)$.

{\rm (\romannumeral2)}~Suppose we are in case {\rm (c)}. Then $\Hg(X) =
\Sp_D(V,\phi)$. The Hodge ring $\cB^\gdot(X)$ is generated by the divisor
classes together with the spaces of Weil classes $W_k \subset \cB^2(X)$, where
$k$ runs through the set of imaginary quadratic fields contained in $D$.

{\rm (\romannumeral3)}~Suppose we are in case {\rm (d)}. Then the Hodge ring
$\cB^\gdot(X)$ is generated by divisor classes, i.e., $\cB^\gdot(X) =
\cD^\gdot(X)$. Either $\Hg(X) = \Sp(V,\phi)$, in which case $\cB^\gdot(X^n) =
\cB^\gdot(X^n)$ for all $n$, or $\Hg(X)$ is isogenous to a $\mQ$-form of $\SL_2
\times \SL_2 \times \SL_2$, in which case there are exceptional Hodge classes
in $\cB^2(X^2)$. In the latter case these exceptional Hodge classes are not of
Weil type.

{\rm (\romannumeral4)}~Suppose we are not in one of the cases {\rm (a)}, {\rm
(b)}, {\rm (c)} or {\rm (d)}. Then $\Hg(X) = \Sp_D(V,\phi)$ and $\cB^\gdot(X^n)
= \cD^\gdot(X^n)$ for all $n$.}
\whiteB

Let us note that in the cases (a), (b) and (c) the Weil classes are really
needed to generate the Hodge ring; in these cases we have $\cD^2(X) \neq
\cB^2(X)$. See [\refMZCrel], especially Example 8 and Criterion 13.
\whiteC

Next we consider some special cases in dimension 5.
\whiteC

(e)~The abelian variety $X$ is isogenous to a product $X_1^2 \times X_2$, where
$X_1$ and $X_2$ are as in (a).
\whiteC

(f)~The abelian variety $X$ is isogenous to a product $X_0 \times X_1 \times
X_2$, where $X_0$ is an elliptic curve, where $X_1$ and $X_2$ are as in (a),
and such that $X_0$ and $X_1$ are not isogenous.
\whiteC

(g)~The abelian variety $X$ is isogenous to a product $X_1 \times X_2$ where
$X_1$ is an elliptic curve with complex multiplication by an imaginary
quadratic field $k$ and where $X_2$ is a simple abelian fourfold such that
there exists an embedding $k \hookrightarrow \Eo(X_2)$ via which $k$ acts on
$T_{X_2,0}$ with multiplicities $(1,3)$.
\whiteB

\advance\secno by 1
\edef\refThmTwo{\number\chno.\number\secno}
\noindent
{\bf (\chapno.\sectno) Theorem.}\enspace {\sl Let $X$ be a complex abelian
variety of dimension 5. Let $V$, $\phi$, $D$ and $\Sp_D(V,\phi)$ have the same
meaning as in {\rm (\refMainThm)}.

{\rm (\romannumeral1)}~Suppose we are in case {\rm (e)}. Consider the space of
Weil classes $W_k \subset \cB^2(X_1 \times X_2)$ and write $W_{k,\alpha}
\subset \cB^2(X)$ for its image under the map $\cB^2(X_1 \times X_2)
\rightarrow \cB^2(X)$ induced by a surjective homomorphism $\alpha \colon X
\rightarrow X_1 \times X_2$. Then the Hodge ring $\cB^\gdot(X)$ is generated by
the subalgebra $\cD^\gdot(X)$ of divisor classes together with the subspaces
$W_{k,\alpha}$. The Hodge group $\Hg(X)$ is strictly contained in
$\Sp_D(V,\phi)$.

{\rm (\romannumeral2)}~Suppose we are in case {\rm (f)}. Then $\Hg(X) =
\Hg(X_0) \times \Hg(X_1 \times X_2)$. For every $n \geq 1$ the Hodge ring
$\cB^\gdot(X^n)$ is generated by the images of $\cB^\gdot(X_0^n)$ and
$\cB^\gdot(X_1^n \times X_2^n)$. In particular, $\cB^\gdot(X)$ is generated by
the divisor classes $\cD^\gdot(X)$ together with the pull-backs of the Weil
classes in $W_k \subset \cB^2(X_1 \times X_2)$.

{\rm (\romannumeral3)}~Suppose we are in case {\rm (g)}. Then the Hodge ring
$\cB^\gdot(X)$ is generated by divisor classes, i.e., $\cB^\gdot(X) =
\cD^\gdot(X)$. The Hodge group $\Hg(X)$ is strictly contained in
$\Sp_D(V,\phi)$.

{\rm (\romannumeral4)}~Suppose we are not in one of the cases {\rm (e)}, {\rm
(f)} or {\rm (g)}. Decompose $X$, up to isogeny, as a product of elementary
abelian varieties, say $X \sim Y_1^{m_1} \times \cdots \times Y_r^{m_r}$. Then
$\Hg(X) = \Hg(Y_1^{m_1}) \times \cdots \Hg(Y_r^{m_r})$. For every $n \geq 1$
the Hodge ring $\cB^\gdot(X^n)$ is generated by the images of the Hodge rings
$\cB^\gdot(Y_j^{m_j})$. In particular, if $X$ has no simple factor of dimension
4 then $\Hg(X) = \Sp_D(V,\phi)$ and $\cB^\gdot(X^n) = \cD^\gdot(X^n)$ for every
$n \geq 1$.}
\whiteB

In the decomposition (up to isogeny) $X \sim Y_1^{m_1} \times \cdots \times
Y_r^{m_r}$ in (\romannumeral4) we require the $Y_j$ to be simple, pairwise
non-isogenous, and the $m_j$ are positive integers. Further we remark that in
the cases (e) and (f) the pull-backs of the Weil classes are needed to
generate the Hodge ring of $X$; in these cases we have $\cD^2(X) \neq \cB^2(X)$
and $\cD^3(X) \neq \cB^3(X)$.
\whiteC

As pointed out at the beginning of the introduction, the above results were
already known for {\sl simple\/} abelian varieties. In the
present paper we are therefore mainly concerned with non-simple abelian
varieties. We prove some lemmas which in certain cases allow us to determine
the Hodge group of a product $X_1 \times X_2$, knowing the Hodge groups
$\Hg(X_i)$ of the factors. Using these results we shall determine the Hodge
groups of all complex abelian varieties $X$ with $\dim(X) \leq 5$.
\whiteC

The paper is organised as follows. In the first section we review the notion of
a Hodge group and we recall a number of properties that we shall use. In \S 2
we give an overview of the situation for simple abelian varieties of low
dimension. In \S 3 we prove a couple of general lemmas which allow us to
analyse certain product situations. In \S 4 we analyse Hodge groups of simple
abelian surfaces of CM-type. Putting everything together the main theorems are
proven in \S 5.
\whiteB\whiteB

\advance\chno by 1
\secno=0
\noindent
\centerline{{\bf \S \chapno.~Hodge groups of abelian varieties.}}
\whiteA

\noindent
\advance\secno by 1
{\bf (\chapno.\sectno)}\enspace Let $X$ be an abelian variety over an
algebraically closed field $k$. Set $D = \End^0(X) := \End(X) \otimes_\mZ \mQ$.
A polarization of $X$ induces a positive (Rosati-) involution, say $d \mapsto
d^\dagger$, of $D$.

Now assume that $X$ is simple. Then $D$ is a division algebra and we have $D
\supset F \supset F_0 \supset \mQ$ with
$$
F = \Cent(D)\, ,\qquad F_0 =\{a \in F\mid a^\dagger = a\}\, .
$$
We write
$$
e_0 = [F_0:\mQ]\, ,\qquad e=[F:\mQ]\, ,\qquad d^2 = [D:F]\, .
$$
By the classification due to Albert (see [\refMAV], \S~21) the division algebra
$D$ is of one of the following types.
\whiteC

\setbox0=\hbox{Type \Romno 4($e_0,d$):\ \ }
\noindent
\hangindent=\wd0
\hangafter=1
\hbox to\wd0{Type \Romno 1($e_0$):\hfil}$e=e_0$, $d=1$; $D = F = F_0$ is a
totally real field.
\smallskip

\noindent
\hangindent=\wd0
\hangafter=1
\hbox to\wd0{Type \Romno 2($e_0$):\hfil}$e=e_0$, $d=2$; $D$ is a quaternion
algebra over a totally real field $F=F_0$; $D$ splits at all infinite places.
\smallskip

\noindent
\hangindent=\wd0
\hangafter=1
\hbox to\wd0{Type \Romno 3($e_0$):\hfil}$e=e_0$, $d=2$; $D$ is a quaternion
algebra over a totally real field $F=F_0$; $D$ is inert at all infinite places.
\smallskip

\noindent
\hangindent=\wd0
\hangafter=1
\hbox to\wd0{Type \Romno 4($e_0,d$):\hfil}$e=2e_0$; $F$ is a CM-field with
totally real subfield $F_0$; $D$ is a division algebra of rank $d^2$ over $F$.
\whiteC

We say that a (simple) abelian variety $X$ is of Type $A$ (with $A \in
\{{\rm \Romno 1, \Romno 2,\Romno 3, \Romno 4}\}$) if $\End^0(X)$ is an algebra
of the corresponding type.

We refer to [\refOortEnd] for results about which algebras in the Albert
classification occur as the endomorphism algebra of an abelian variety. (Note
that there is a misprint in Table 8.1 of [\refOortEnd]; the author informs us
that in the last line of this table it should read: ``occurs if and only if
$2g/ed^2 \geq 1$ but {\sl excluded\/} \Romno 4($1,1$), $g=2$ and \Romno
4($1,1$), $g=4$.'')
\whiteB

\noindent
\advance\secno by 1
{\bf (\chapno.\sectno)}\enspace Let $X$ be a complex abelian variety, $X \neq
0$. We write $V = V_X = H_1(X(\mC),\mQ)$, which is a polarizable $\mQ$-Hodge
structure of type $(-1,0) + (0,-1)$. This Hodge structure can be described by
giving a homomorphism of algebraic groups over $\mR$
$$
h \colon \mS \rightarrow \GL(V)_\mR\, ,
$$
where $\mS = \Res_{\mC/\mR} \mG_{m,\mC}$.

The Mumford-Tate group $\MT(X)$ of $X$ is defined to be the smallest algebraic
subgroup $M \subset \GL(V)$ (over $\mQ$) such that $h$ factors through $M_\mR$.
In practice it is often more convenient to work with the {\sl Hodge group\/}
$\Hg(X)$. We can define it by $\Hg(X) = \MT(X) \cap \SL(V)$. For a more direct
definition, consider the $\mR$-subtorus $U^1 \subset \mS$ given on points by
$$
U^1(\mR) = \{z \in \mC^\ast \mid z \bar{z} = 1\} \subset \mC^\ast = \mS(\mR)\,
{}.
$$
Then $\Hg(X)$ is the smallest algebraic subgroup $H \subset \GL(V)$ such that
the restriction of $h$ to $U^1$ factors through $H_\mR$.

The Mumford-Tate group $\MT(X)$ contains the torus $\mG_{m,\mQ} \subset \GL(V)$
of homotheties. The group $\MT(X)$ is the almost direct product of
$\mG_{m,\mQ}$ and $\Hg(X)$.

The Hodge group $\Hg(X)$ is a connected reductive algebraic group. Viewing $D =
\Eo(X)$ as a subalgebra of $\End_\mQ(V)$ we have $D = \End_\mQ(V)^{\Hg(X)}$. If
$\phi \colon V \times V \rightarrow \mQ$ is the Riemann form associated to a
polarization of $X$ (so $\phi$ is a symplectic form) then
$$
\Hg(X) \subset \Sp_D(V,\phi)\, ,
$$
the centralizer of $D$ in the symplectic group $\Sp(V,\phi)$.

The Hodge group $\Hg(X)$ is a torus if and only if $X$ is of CM-type. If $X$
has no factors of Type \Romno 4 then $\Hg(X)$ is semi-simple. (See
[\refMumNote], \S 2 and [\refTank], Lemma 1.4.)

For $n \geq 1$ we can identify $\Hg(X^n)$ with $\Hg(X)$, acting diagonally on
$V_{X^n} = (V_X)^n$. More generally, if $n_1, \ldots, n_r \in \mZ_{\geq 1}$
then we can identify $\Hg(X_1^{n_1} \times \cdots \times X_r^{n_r})$ with
$\Hg(X_1 \times \cdots \times X_r)$.
\whiteB

\noindent
\advance\secno by 1
\edef\refHgLie{\number\chno.\number\secno}%
{\bf (\chapno.\sectno)}\enspace Write $\hg(X)$ for the Lie algebra of $\Hg(X)$.
If $W$ is a $\Hg(X)$-module then $W^{\Hg(X)} = W^{\hg(X)}$, since $\Hg(X)$ is
connected. Thus, for instance, $\Eo(X)$ can be computed as the
$\hg(X)$-invariants in $\End_\mQ(V)$.

The following description of $\hg(X)$ proves to be very useful. We have a Hodge
decomposition $V_\mC = V_\mC^{-1,0} \oplus V_\mC^{0,-1}$. Let the endomorphism
$J = J_X \in \End(V_\mC)$ be given by
$$
J_X(v) =
\cases{iv,&if $v \in V_\mC^{-1,0}$,\cr
-iv,&if $v \in V_\mC^{0,-1}$.\cr}
$$
Note that $J_X^2 = -\id$. Then $\hg(X) \subset \End(V)$ is the smallest
$\mQ$-Lie
subalgebra $\gh \subset
\End(V)$ such that $\gh_\mC$ contains $J_X$; see [\refZarLNM]. In fact, since
$V_\mC^{-1,0}$ and $V_\mC^{0,-1}$ are complex conjugate we even have $J_X \in
\hg(X)_\mR$.

We remark that the same automorphism $J_X$ can also be viewed as the element
$h(i) \in \Hg(X)(\mR)$. (This element is usually referred to as the Weil
operator.)
\whiteB

\noindent
\advance\secno by 1
{\bf (\chapno.\sectno)}\enspace The cohomology ring $H^\gdot(X,\mQ)$ is
naturally isomorphic to the exterior algebra on $V^\vee$. The Hodge group
$\Hg(X)$ acts on this ring. The $\Hg(X)$-invariants in $H^\gdot(X,\mQ)$ are
precisely the Hodge classes. Writing $\cB^i(X) \subset H^{2i}(X,\mQ)$ for the
subspace of Hodge classes we obtain a graded $\mQ$-algebra $\cB^\gdot(X) =
\oplus_i \cB^i(X)$, called the Hodge ring of $X$.

The Hodge classes in $H^2(X,\mQ)$ (i.e., the elements of $\cB^1(X) =
H^2(X,\mQ)^{\Hg(X)}$) are called the divisor classes. We write $\cD^\gdot(X)
\subset \cB^\gdot(X)$ for the $\mQ$-subalgebra generated by the divisor
classes. The Hodge classes in $\cD^\gdot(X)$ are called the {\sl
decomposable\/} Hodge classes. The elements of $\cB^\gdot(X)$ not in
$\cD^\gdot(X)$ are called {\sl exceptional\/} Hodge classes.
\whiteB

\noindent
\advance\secno by 1
\edef\refMinWtsA{\number\chno.\number\secno}%
{\bf (\chapno.\sectno)}\enspace Consider the tautological representation $\rho
\colon \hg(X) \rightarrow \End(V_X)$. The fact that $V_X$ is a polarizable
Hodge structure of weight 1 puts strong restrictions on this representation. We
shall summarize this here; for further details we refer to [\refDelShim], \S 1.
See also [\refPink], \S 4 and [\refZarWts].

Consider the decomposition
$$
\hg(X) \otimes \mR = \gc \times \gg_1 \times \cdots \times \gg_q
$$
of $\hg(X) \otimes \mR$ as a product of its center $\gc$ and a number of
$\mR$-simple factors $\gg_i$. A certain number of these factors, say $\gg_1,
\ldots,\gg_r$ are non-compact. (Here $0 \leq r \leq q$.) As remarked above,
$J_X$ can also be viewed as the Weil operator in $\Hg(X)(\mR)$. The
automorphism $\Ad(J_X)$ of $\Hg(X)(\mR)$ is a Cartan involution (see
[\refDelWeilK3], \S 2, especially Lemma 2.8 and Proposition 2.11). This implies
that each $\mR$-simple factor $\gg_i$ is a form of a compact real Lie algebra
and is therefore absolutely simple.

Now consider the representation $\rho$. Let $W \subset V_X \otimes \mC$ be an
irreducible $\hg(X) \otimes \mC$-submodule. Then $W$ decomposes as an external
tensor product
$$
W = \chi_0 \boxtimes W_1 \boxtimes \cdots \boxtimes W_q\, ,
$$
where $\chi_0$ is a character of $\gc$ and where $W_i$ is an irreducible
representation of $\gg_i \otimes_\mR \mC$. With these notations, the Lie
algebra $\hg(X)$ and its representation $\rho$ have the property that

(\romannumeral1)~all simple factors $\gg_i$ are of classical type A$_\ell$,
B$_\ell$, C$_\ell$ or D$_\ell$,

\noindent
and for every irreducible $\hg(X) \otimes \mC$-submodule $W \subset V_\mC$ as
above, we have

(\romannumeral2)~at most one of the representations $W_1, \ldots, W_r$ is
non-trivial (with $r$ as introduced above),

(\romannumeral3)~if $W_i$ is a non-trivial $\gg_i$-module ($1 \leq i \leq q$)
then its highest weight (w.r.t.\ a chosen Cartan subalgebra of $\gg_i$ and a
choice of a basis for the root system) is miniscule in the sense of
[\refBourLie], Chap.~8, \S 7, n$^{{\rm o}}$ 3.

These facts can be found in [\refDelShim], sections 1.3 and 2.3. For the
purpose of this paper it actually suffices to know that in every irreducible
$\hg(X) \otimes \mC$-module $W$ {\sl at most\/} one of the non-compact factors
$W_1, \ldots W_r$ is non-trivial. Let us sketch the argument. Decompose $\hg(X)
\otimes \mR = \gc \times \gg_1 \times \cdots \times \gg_q$ as above and write
$J = (J_0,J_1,\ldots,J_q)$. The fact that $\Ad(J_X)$ is a Cartan involution
implies that $J_1, \ldots,J_r$ are all non-zero. If $W_i$ ($1 \leq i \leq q$)
is a non-trivial $\gg_i$-module and $J_i \neq 0$ then the simplicity of $\gg_i$
implies that $J_i$ has trace 0 on $W_i$ and has therefore at least 2 different
eigenvalues. Combining everything we see that if there are two non-compact
factors, say $\gg_i$ and $\gg_j$ ($1 \leq i < j \leq r$), acting non-trivially
on $W$ then the operator $J$ has at least 3 different eigenvalues, as
$(J_i,J_j)$ has at least 3 different eigenvalues on $W_i \boxtimes W_j$.
Contradiction.

Let us remark that if $X$ is of CM-type then $\Hg(X)$ is a torus and we have $q
= 0$ in the above. (Thus, statements (\romannumeral1), (\romannumeral2) and
(\romannumeral3) become void in this case.) Next suppose that $X$ itself is not
of CM-type but that it contains an abelian subvariety $Y$ of CM-type. Then the
semi-simple part of $\hg(X)$ acts trivially on $V_Y$; in particular we find
$\hg(X) \otimes \mC$-modules $W$ as above for which all factors $W_i$ ($1 \leq
i\leq q$) are trivial. On the other hand, if $X$ does not contain an abelian
subvariety of CM-type then it can be shown that for every $W$ as above
precisely one of the factors $W_1, \ldots W_r$ is non-trivial.
\whiteB

\noindent
\advance\secno by 1
\edef\refMinWtsB{\number\chno.\number\secno}%
{\bf (\chapno.\sectno)}\enspace Let $\gh$ be a reductive Lie algebra over
$\mQ$. We shall say that $\gh$ is of non-compact type if $\gh \otimes \mR$ does
not have compact simple factors.

Suppose $\gh$ is of non-compact type. Let $\rho \colon
\gh \rightarrow \End(V)$ be a finite-dimensional representation. The Lie
algebra $\gh \otimes \mC$ decomposes as $\gh_\mC = \gc \oplus \gg_1 \oplus
\cdots \oplus \gg_q$, where $\gc$ is its center and where the $\gg_i$ are its
simple factors. As in (\refMinWtsA), every irreducible $\gh_\mC$-submodule $W
\subset V_\mC$ decomposes as an external tensor product $W = \chi_0 \boxtimes
W_1 \boxtimes \cdots \boxtimes W_q$. We shall say that $\rho$ is a {\sl length
1 representation of non-compact type\/} if all simple factors $\gg_i$ are of
classical type and if every irreducible $\gh \otimes \mC$-submodule $W \subset
V_\mC$ satisfies the conditions

(\romannumeral2$^\prime$)~at most one of the representations $W_1, \ldots, W_q$
is non-trivial,

(\romannumeral3)~if $W_i$ is a non-trivial $\gg_i$-module then its highest
weight is miniscule.

Our terminology is based on [\refZarWts], where the length of an irreducible
representation of a simple Lie algebra is defined. As remarked in loc.\ cit.,
2.1, such a representation has length 1 precisely if the Lie algebra is of
classical type and the highest weight of the representation is miniscule. (See
also [\refSerGAHT], \S 3.)

Needless to say, our interest in length 1 representations comes from the facts
recalled in (\refMinWtsA). If $X$ is an abelian variety such that $\hg(X)$ is
of non-compact type then these facts tell us that the tautological
representation $\rho \colon \hg(X) \rightarrow \End(V_X)$ is a length 1
representation of non-compact type. More generally, if the Hodge Lie algebra
$\hg(X)$ contains an ideal $\gg$ which is of non-compact type then the
restricted representation $\rho_{|\gg}$ is a length 1 representation of
non-compact type.
\whiteB

\noindent
\advance\secno by 1
\edef\refUniqRep{\number\chno.\number\secno}%
{\bf (\chapno.\sectno) Remark.}\enspace Later in the paper we shall consider
$\mQ$-Lie algebras $\gh$ of non-compact type for which there is a unique
faithful irreducible  length 1 representation of non-compact type (up to
isomorphism). For instance, let $\gh$
be a simple $\mQ$-Lie algebra of non-compact type. Then there exists a number
field $K$ and an absolutely simple $K$-Lie algebra $\gg$ such that $\gh \cong
\Res_{K/\mQ}\,
\gg$. Writing $\Sigma_K$ for the set of embeddings of $K$ into $\mC$ we have
$\gh_\mC = \oplus_{\sigma \in \Sigma_K} \gg_{(\sigma)}$, where $\gg_{(\sigma)}
= \gg \otimes_{K,\sigma} \mC$. We claim that if the (absolute) root system of
$\gg$ is of type C$_\ell$ ($\ell \geq 1$) then $\gh$ has a unique irreducible
representation of length 1.

To see this, let us first remark that a simple Lie algebra of type C$_\ell$
($\ell \geq 1$) over $\mC$ has a unique irreducible representation with
miniscule highest weight, see [\refBourLie], Chap.~8, \S 7, n$^{{\rm o}}$ 3.
Now write $\Sigma_K = \{\sigma_1, \ldots, \sigma_r\}$ and let $V_{(i)}$ ($1
\leq i \leq r$) be the irreducible $\gh_\mC$-module which is irreducible as a
$\gg_{(\sigma_i)}$-module with miniscule highest weight, and on which the
factors $\gg_{(\sigma_j)}$ with $i \neq j$ act trivially. If $\rho \colon \gh
\rightarrow \End(V)$ is an irreducible length 1 representation of non-compact
type then
$$
V_\mC \cong V_{(1)}^{m_1} \oplus \cdots \oplus V_{(r)}^{m_r}
$$
for certain multiplicities $m_i$. But if $L$ is the normal closure of $K$
inside $\mC$ then $\Gal(L/\mQ)$ permutes the factors $\gg_{(\sigma_i)}$
transitively, and it follows from the fact that $V_\mC$ is defined over $\mQ$
that we must have $m_1 = m_2 = \cdots = m_r$. Therefore, if $\rho^\prime \colon
\gh \rightarrow \End(V^\prime)$ is another irreducible length
1 representation of non-compact type then there is a relation $(\rho_\mC)^M
\cong (\rho^\prime_\mC)^N$ for certain
integers $M$ and $N$. But this is possible only if $\rho \cong \rho^\prime$.

That, conversely, every $\gh$ of non-compact type as above has an irreducible
(symplectic) length 1 representation of non-compact type can be seen from the
description of such $\gh$'s in
terms of algebras with involution, as in [\refJacobs], Chap.\ \Romno{10}.
\whiteB

\noindent
\advance\secno by 1
\edef\refHazMur{\number\chno.\number\secno}%
{\bf (\chapno.\sectno)}\enspace Consider the following condition on the complex
abelian variety $X$:
$$
\cB^\gdot(X^n) = \cD^\gdot(X^n)\qquad \hbox{{\rm for all $n$}}\, .\leqno{\rm
(D)}
$$
If this condition is satisfied then the Hodge conjecture is ``trivially'' true
for all $X^n$.

As was recalled above, the Hodge group $\Hg(X)$ is contained in the algebraic
group
$\Sp_D(V,\phi)$. It was shown by Hazama [\refHazHC] and Murty [\refMurtyHC]
(independently) that
$$
\Hg(X) = \Sp_D(V,\phi)
\quad\Longleftrightarrow\quad
\left(
\vcenter{
\setbox0=\hbox{$X$ has no factors of type \Romno 3}
\hbox{$X$ has no factors of type \Romno 3}
\hbox to \wd0{\hfil{\sl and\/}\hfil}
\hbox to \wd0{\hfil $\cD^\gdot(X^n) = \cB^\gdot(X^n)$ for all $n$\hfil}
}
\right)\, .
$$
\whiteB

\noindent
\advance\secno by 1
{\bf (\chapno.\sectno)}\enspace Let $K$ be a subfield of $\Eo(X)$, with $1 \in
K$ acting as the identity on $X$. Write $\Sigma_K$ for the set of embeddings of
$K$ into $\mC$. Let $T_{X,0}$ be the tangent space of $X$ at the origin. The
action of (an order of) $K$ on $X$ makes $T_{X,0}$ into a module under $K
\otimes_\mQ \mC = \prod_{\sigma \in \Sigma_K} \mC$. This gives a decomposition
$$
T_{X,0} = \bigoplus_{\sigma \in \Sigma_K} T^{(\sigma)}\, .
$$
Let $n_\sigma = \dim_\mC T^{(\sigma)}$. If $\overline{\sigma} \colon K
\rightarrow \mC$ is the complex conjugate of $\sigma$ then $n_\sigma +
n_{\overline{\sigma}} = r := 2 \dim(X)/[K:\mQ]$.

If $K$ is imaginary quadratic then we say that it acts on $T_{X,0}$ with
multiplicities $(a,b)$ if $n_\sigma = a$, $n_{\overline{\sigma}} = b$ for some
ordering $\Sigma_K = \{\sigma,\overline{\sigma}\}$.

The inclusion $K \subset \Eo(X)$ induces on $V_X$ the structure of an
$r$-dimensional $K$-vector space. The 1-dimensional $K$-vector space $W_K =
W_K(X) := \wedge^r_K V_X^\vee$ can be identified in a natural way with a
subspace of $H^r(X,\mQ)$; we call $W_K$ the space of Weil classes w.r.t.\ $K$.
(We refer the reader to [\refWeil].) It is known that either $W_K$ consists
entirely of Hodge classes or $0 \in W_K$
is the only Hodge class in $W_K$. Whether $W_K$ consists of Hodge classes and,
if so, whether these classes are exceptional or not, can be answered purely in
terms of the data $K \subset \Eo(X)$ and the action of $K$ on $T_{X,0}$, see
[\refMZCrel]. For instance, it is shown there that $W_K$ consists of Hodge
classes if and only if $n_\sigma = n_{\overline{\sigma}}$ for all $\sigma \in
\Sigma_K$. Note also that the Hodge Lie algebra is contained in the Lie algebra
$\End_K(V_X)$ of $K$-linear endomorphisms of $V_X$ and that it acts on $W_K$
through the $K$-linear trace map $\trace_K \colon \End_K(V_X) \rightarrow K$.
In particular, $W_K$ consists of Hodge classes precisely if $\hg(X) \subseteq
\gsl_K(V_X)$.

For later use, let us note the following. Suppose $X$ is isogenous to a
product, say $X \sim X_1 \times X_2$. Then (an order of) $K$ acts on both $X_1$
and $X_2$. Let $r_1$ ($i=1,2$) be the $K$-dimension of $V_{X_i}$, so that $r =
r_1 + r_2$. We have associated spaces of Weil classes $W_{K}(X_1) \subset
H^{r_1}(X_1,\mQ)$ and $W_K(X_2) \subset H^{r_2}(X_2,\mQ)$. Viewing
$H^{r_1}(X_1,\mQ) \otimes H^{r_2}(X_2,\mQ)$ as a subspace of $H^r(X,\mQ)$ via
the K\"unneth decomposition, the space of Weil classes $W_K(X)$ can naturally
be identified with $W_K(X_1) \otimes_K W_K(X_2)$; see also [\refMZCrel],
section 7.
\whiteB\whiteB

\advance\chno by 1
\secno=0
\noindent
\centerline{{\bf \S \chapno.~Simple abelian varieties of dimension $\leq 5$.}}
\whiteA

\noindent
We shall give a short overview of the situation for simple complex abelian
varieties of low dimension. Thus, in this section we shall assume $X$ to be
{\sl simple}.

For $g := \dim(X) \leq 3$ and $g = 5$ we always find that $\Hg(X) =
\Sp_{D}(V,\phi)$. Since type \Romno 3 does not occur for $g \leq 3$ and $g = 5$
($X$ simple!), it follows that $\cB^\gdot(X^n) = \cD^\gdot(X^n)$ for all $n$.
(See (\refHazMur).) In particular the Hodge conjecture is true for all such
$X^n$. A useful references for the results stated below is [\refRibet].

We shall give an overview of the cases that occur. If $F$ is a CM-field with
totally real subfield $F_0$ and complex conjugation $x \mapsto \bar{x}$ then we
shall write $\UU_F$ for the algebraic torus over $\mQ$ given on points by
$$
\UU_F(R) = \{x \in (F \otimes_\mQ R)^\ast \mid x \bar{x} = 1\}\, .
$$
\whiteB

\noindent
\advance\secno by 1
\edef\refEllC{\number\chno.\number\secno}%
{\bf (\chapno.\sectno) ${\bfit g}$=1.}\enspace There are two cases to
distinguish.
\whiteC

\setbox0=\hbox{Type \Romno 4(3,1):\ \ }
\noindent
\hangindent=\wd0
\hangafter=1
\hbox to\wd0{Type \Romno 1(1):\hfil}$X$ is an elliptic curve with $\Eo(X) =
\mQ$. Then $\Hg(X) = \Sp(V,\phi) \cong \SL_{2,\mQ}$.
\whiteC

\noindent
\hangindent=\wd0
\hangafter=1
\hbox to\wd0{Type \Romno 4(1,1):\hfil}$X$ is an elliptic curve with CM by an
imaginary quadratic field $F$. Then $\Hg(X) = \UU_F$.
\whiteB

\noindent
\advance\secno by 1
{\bf (\chapno.\sectno) ${\bfit g}$=2.}\enspace There are four cases.
\whiteC

\noindent
\hangindent=\wd0
\hangafter=1
\hbox to\wd0{Type \Romno 1(1):\hfil}$X$ is an abelian surface with $\Eo(X) =
\mQ$. Then $\Hg(X) = \Sp(V,\phi) \cong \Sp_{4,\mQ}$.
\whiteC

\noindent
\hangindent=\wd0
\hangafter=1
\hbox to\wd0{Type \Romno 1(2):\hfil}$\Eo(X) = F$ is a real quadratic field.
Then there is a unique $F$-symplectic form $\psi \colon V \times V \rightarrow
F$ such that $\phi = \trace_{F/\mQ} \psi$. The Hodge group is given by $\Hg(X)
= \Res_{F/\mQ} \Sp_{F}(V,\psi)$.
\whiteC

\noindent
\hangindent=\wd0
\hangafter=1
\hbox to\wd0{Type \Romno 2(1):\hfil}$D = \Eo(X)$ is a quaternion algebra over
$\mQ$, split at $\infty$. Write $D^\opp$ for the opposite algebra, and let $x
\mapsto x^\ast$ be the canonical involution. Then $\Hg(X)$ is the algebraic
group $\UU_{D^\opp}$ given on points by $\UU_{D^\opp}(\mQ) = \{x \in
(D^\opp)^\ast \mid x x^\ast = 1\}$.
\whiteC

\noindent
\hangindent=\wd0
\hangafter=1
\hbox to\wd0{Type \Romno 4(2,1):\hfil}$\Eo(X) = F$ is a quartic CM-field not
containing an imaginary quadratic subfield. We have $\Hg(X) = \UU_F$.
\whiteB

\noindent
\advance\secno by 1
\edef\refThreef{\number\chno.\number\secno}%
{\bf (\chapno.\sectno) ${\bfit g}$=3.}\enspace There are four cases.
\whiteC

\noindent
\hangindent=\wd0
\hangafter=1
\hbox to\wd0{Type \Romno 1(1):\hfil}$X$ is an abelian 3-fold with $\Eo(X) =
\mQ$. Then $\Hg(X) = \Sp(V,\phi) \cong \Sp_{6,\mQ}$.
\whiteC

\noindent
\hangindent=\wd0
\hangafter=1
\hbox to\wd0{Type \Romno 1(3):\hfil}$\Eo(X) = F$ is a totally real cubic field.
There is a unique $F$-symplectic form $\psi \colon V \times V \rightarrow F$
such that $\phi = \trace_{F/\mQ} \psi$. The Hodge group is given by $\Hg(X) =
\Res_{F/\mQ} \Sp_{F}(V,\psi)$.
\whiteC

\noindent
\hangindent=\wd0
\hangafter=1
\hbox to\wd0{Type \Romno 4(1,1):\hfil}$\Eo(X) = F$ is an imaginary quadratic
field; given $a \in F$ with $\bar{a} = -a$ there is a unique $F$-hermitian form
$\psi \colon V \times V \rightarrow F$ such that $\phi = \trace_{F/\mQ}(a \cdot
\psi)$ and $\Hg(X) = \UU_F(V,\psi)$.
\whiteC

\noindent
\hangindent=\wd0
\hangafter=1
\hbox to\wd0{Type \Romno 4(3,1):\hfil}$\Eo(X) = F$ is a CM-field of degree 6
over $\mQ$. Then $\Hg(X) = \UU_F$.
\whiteB

\noindent
\advance\secno by 1
\edef\refPropLowg{\number\chno.\number\secno}%
{\bf (\chapno.\sectno) Proposition.}\enspace {\sl Let $X$ be a simple complex
abelian variety with $g = \dim(X) \leq 3$.

{\rm (\romannumeral1)}~The Hodge Lie algebra $\hg(X)$ is of non-compact type in
the sense of {\rm (\refMinWtsB)}.

{\rm (\romannumeral2)}~Suppose $X$ is of CM-type. Then $\Hg(X)$ is a
$g$-dimensional algebraic torus. It is $\mQ$-simple, except when $\dim(X) = 3$
and the sextic CM-field $\Eo(X)$ contains an imaginary quadratic field.

{\rm (\romannumeral3)}~Suppose $X$ is not of CM-type. Then $\Hg(X)$ is a
$\mQ$-simple algebraic group, except when $\dim(X) = 3$ and $\Eo(X)$ is an
imaginary quadratic field. If $\Hg(X)$ is $\mQ$-simple then (up to isomorphism)
there is exactly one faithful irreducible representation of $\hg(X)$ over $\mQ$
which is of length 1.}
\whiteD

\Proof~Most of the claims are easily read off from the above. For
(\romannumeral1) let us add that if $g=3$ and $\Eo(X) = F$ is imaginary
quadratic (Type \Romno 4(1,1)), $F$ necessarily acts on the tangent space with
multiplicities (2,1). (An action with multiplicities $(3,0)$ is excluded; see
[\refShimFam], Proposition~14.) Thus $\Hg(X)_\mR$ is a unitary group of
signature $(2,1)$, which has a non-compact $\mR$-simple derived group. For
(\romannumeral2), use Lemma (3.7) below. For the last assertion of
(\romannumeral3) one uses (\refUniqRep). $\square$
\whiteB

\noindent
\advance\secno by 1
{\bf (\chapno.\sectno) ${\bfit g}$=4.}\enspace The case $g =4$ is more involved
and was studied in [\refMZDuke]. In particular, in op.\ cit.\ we already proved
Theorem~(\refMainThm) for simple abelian fourfolds. (This covers the cases (b),
(c) and (d) of the introduction.) We here only recall some of the most
interesting cases.

(\romannumeral1)~For $g=4$ it is no longer true that $\Hg(X)$ is determined by
$\Eo(X)$ together with its action on the tangent space at the origin. Namely,
if $g=4$ and $\Eo(X) = \mQ$ then either $\Hg(X) = \Sp(V,\phi) \cong
\Sp_{8,\mQ}$, or $\Hg(X)$ is a $\mQ$-form of an almost direct product of three
copies of $\SL_2$. (See [\refMumNote].) In both cases the Hodge ring of $X$ is
generated by divisor classes, but if $\Hg(X)$ is isogenous to a $\mQ$-form of
$\SL_2^3$ then there are exceptional Hodge classes in $H^4(X^2,\mQ)$.

(\romannumeral2)~For $g = 4$ we find cases where in addition to divisor classes
we also need Weil classes to generate the Hodge ring. This happens if $\Eo(X)$
contains an imaginary quadratic field $k$ which acts on the tangent space with
multiplicities $(2,2)$. If $X$ is of Type \Romno 3 then this is the case (e.g.,
see [\refMurtyHC], [\refMZCrel]); further it can occur only for $X$ of Type
\Romno 4(1,1) or of Type \Romno 4(4,1). Only in very special cases these Weil
classes are known to be algebraic, see [\refSchoenHC] and [\refvGeemen].
\whiteB

\noindent
\advance\secno by 1
{\bf (\chapno.\sectno) ${\bfit g}$=5.}\enspace As already stated above, $\Hg(X)
= \Sp_{D}(V,\phi)$ for all simple abelian 5-folds. The point here is that 5 is
a prime number, since in fact we have the following result, due to Tankeev
[\refTankPrime]. (See also Ribet's paper [\refRibet].)
\whiteB

\noindent
\advance\secno by 1
\edef\refPrimeDim{\number\chno.\number\secno}%
{\bf (\chapno.\sectno) Theorem.}\enspace {\sl Let $X$ be a simple complex
abelian variety such that $\dim(X)$ is a prime number. Then $\Hg(X) =
\Sp_D(V,\phi)$ and $\cB^\gdot(X^n) = \cD^\gdot(X^n)$ for every $n \geq 1$.}
\whiteB

In connection with this result let us note that a simple $X$ of prime dimension
cannot be of Type \Romno 3, so that the result of Hazama and Murty in
(\refHazMur) applies.
\whiteB\whiteB

\advance\chno by 1
\secno=0
\noindent
\centerline{{\bf \S \chapno.~The Hodge group of a product of abelian
varieties.}}\nobreak
\whiteA\nobreak

\noindent
\advance\secno by 1
\edef\refProd{\number\chno.\number\secno}%
{\bf (\chapno.\sectno)}\enspace Let $X_1$ and $X_2$ be complex abelian
varieties. Write $X = X_1 \times X_2$. Then $\Hg(X)$ is an algebraic subgroup
of $\Hg(X_1) \times \Hg(X_2)$. The two projections $\pr_i \colon \Hg(X)
\rightarrow \Hg(X_i)$ are surjective. From this one easily shows that there
exist Lie algebras $\gg_1$, $\gg_2$, $\gg_3$ and an
automorphism $\phi$ of $\gg_3$ such that
$$
\hg(X_1) \cong \gg_1 \oplus \gg_3\, ,\qquad \hg(X_2) \cong \gg_2 \oplus \gg_3\,
,
$$
and
$$
\def\normalbaselines{\baselineskip20pt\lineskip3pt\lineskiplimit3pt}
\matrix{
\hg(X_1 \times X_2)&&\subseteq&&\hg(X_1) \oplus \hg(X_2)\cr
\wr\Big\Vert&&&&\Big\Vert\wr\cr
\gg_1 \oplus \gg_2 \oplus \Gamma_\phi&\subseteq&\gg_1 \oplus \gg_2 \oplus \gg_3
\oplus \gg_3&\cong&(\gg_1 \oplus \gg_3) \oplus (\gg_2 \oplus \gg_3)\rlap{\quad
,}\cr
}
$$
where $\Gamma_\phi \subseteq (\gg_3 \oplus \gg_3)$ is the graph of the
automorphism $\phi$.

We may have that
$$
\Hg(X_1 \times X_2) \neq \Hg(X_1) \times \Hg(X_2)\, .\eqno(1)
$$
(I.e., $\gg_3 \neq 0$ in the above.) This holds if and only if for some
$m$ and $n$ the Hodge ring $\cB^\gdot(X_1^m \times X_2^n)$ is not generated by
the elements coming from $\cB^\gdot(X_1^m)$ and $\cB^\gdot(X_2^n)$.

In certain cases one can show that an inequality (1) can only hold if
$\Hom(X_1,X_2) \neq 0$. For instance, we have the following result of Hazama
[\refHazNons].
\whiteB

\noindent
\advance\secno by 1
\edef\refHazThm{\number\chno.\number\secno}%
{\bf (\chapno.\sectno) Theorem.} {\sl Let $X_1$ and $X_2$ be complex abelian
varieties which both satisfy condition {\rm (D)} in {\rm (\refHazMur)}.

{\rm (\romannumeral1)}~Suppose $X_1$ and $X_2$ contain no factors of Type
\Romno 4. Then $X_1 \times X_2$ again satisfies {\rm (D)}, and either $\Hom(X_1
,X_2) \neq 0$ or $\Hg(X_1 \times X_2) = \Hg(X_1) \times \Hg(X_2)$.

{\rm (\romannumeral2)}~Suppose $X_1$ has no factors of Type \Romno 4 and $X_2$
is of CM-type. Then $X_1 \times X_2$ again satisfies {\rm (D)} and $\Hg(X_1
\times X_2) = \Hg(X_1) \times \Hg(X_2)$.}
\whiteC

The next lemmas are aimed at proving similar conclusions in other cases.
\whiteB

\noindent
\advance\secno by 1
\edef\refLemIdeal{\number\chno.\number\secno}%
{\bf (\chapno.\sectno) Lemma.}\enspace {\sl Let $X$ be a complex abelian
variety. Suppose that $\hg(X)$ is semi-simple of non-compact type and that, up
to isomorphism, $V_X$ is the only irreducible $\hg(X)$-representation which is
a length 1 representation of non-compact type. Let $Y$ be a simple complex
abelian variety such that $\hg(Y)$ splits as $\hg(Y) = \gg \oplus \gh$;
correspondingly we can write $J_Y = J_1 + J_2$ with $J_1 \in
\gg_\mC$ and $J_2 \in \gh_\mC$. Suppose there exists an isomorphism $\hg(X)
\isomarrow \gg$ with $J_X \mapsto J_1$. Then $\gh = 0$ and $Y$ is isogenous to
$X$.}
\whiteD

\Proof~Write $D = \Eo(X)$ and $F = \Cent(D)$; set $e = [F : \mQ]$ and $d^2 =
\dim_F(D)$. We have $D \otimes_\mQ \mC \cong M_d(\mC)_{(1)} \times \cdots
\times M_d(\mC)_{(e)}$. There are irreducible $\hg(X)_\mC$-modules $U_1,
\ldots, U_e$, pairwise non-isomorphic, such that $V_X \otimes_\mQ \mC \cong
U_1^d \oplus \cdots \oplus U_e^d$ as $\hg(X)_\mC$-modules.

As $\hg(X)$ is semi-simple, the $F$-linear trace map $\trace_F \colon \hg(X)
\subset \End_F(V_X) \rightarrow F$ is zero. It follows that $\hg(X)_\mC$ acts
on each of the summands $U_j^d$ through $\gsl(U_j^d)$. In particular, on each
of the summands $U_j^d$ the operator $J_X$ has $+i$ and $-i$ as its eigenvalues
(as it has zero trace and satisfies $J_X^2 = -\id$).

Fix an isomorphism $\phi \colon \hg(X) \isomarrow \gg$ with $J_X \mapsto J_1$.
Note that there are no non-trivial $\gg$-invariants in $V_Y$, as $(V_Y)^{\gg}$
is a $\hg(Y)$-submodule of $V_Y$ and $Y$ is simple. The assumption that $V_X$
is the only length 1 irreducible $\gg$-module of non-compact type therefore
implies that $V_Y
\cong V_X^q$ as $\gg$-modules, for some $q \geq 1$. (See the remarks at the end
of section (\refMinWtsB).) Then $\gh$ acts on $V_Y$
through an embedding $\gh \hookrightarrow \End_\gg(V_Y) = M_q(D)$. Thus
$$
V_{Y,\mC} \cong U_1^{dq} \oplus \cdots \oplus U_e^{dq}\, ,
$$
as $\gg_\mC$-modules and each of the factors $U_j^{dq}$ is stable under
$\gh_\mC$. If $\lambda$ is an eigenvalue of $J_2$ on $U_j^{dq}$ then we find
that both $i + \lambda$ and $-i + \lambda$ occur as eigenvalues of $J_Y$ on
$U_j^{dq} \subseteq V_{Y,\mC}$. By definition of $J_Y$ this is possible only if
$\lambda = 0$. We conclude that $J_2$ acts trivially on each factor $U_j^{dq}$.
Hence $\gh =0$.

The graph $\Gamma_\phi \subset \hg(X) \times \hg(Y)$ is a $\mQ$-Lie subalgebra
such that $\Gamma_{\phi,\mC} \ni J_{X \times Y} = (J_X,J_Y)$. Therefore, $\hg(X
\times Y) = \Gamma_\phi$ and some multiple of $\phi$ corresponds to an isogeny
from $X$ to $Y$. $\square$
\whiteB

\noindent
\advance\secno by 1
\edef\refLemProdSS{\number\chno.\number\secno}%
{\bf (\chapno.\sectno) Lemma.}\enspace {\sl Let $X_1$ and $X_2$ be nonzero
complex abelian varieties. Write $X = X_1 \times X_2$. Assume that $\hg(X_2)$
is a $\mQ$-simple Lie algebra of non-compact type and that, up to isomorphism,
$V_{X_2}$ is the only irreducible $\hg(X_2)$-module which is a length 1
representation of non-compact type. Then either $\Hg(X) = \Hg(X_1) \times
\Hg(X_2)$ or $\Hom(X_2,X_1) \neq 0$.}
\whiteD

\Proof~Assume that $\Hg(X) \neq \Hg(X_1) \times \Hg(X_2)$. Using the notations
of (\refProd) the assumption that $\hg(X_2)$ is $\mQ$-simple implies that
$\hg(X) = \gg_1 \oplus \gg_3 \isomarrow \hg(X_1)$ and $\hg(X_2) \cong \gg_3$.

There exists a {\sl simple\/} abelian subvariety $Y \subset X_1$ such that the
ideal $\gg_3 \subset \hg(X_1)$ acts non-trivially on $V_Y \subset V_{X_1}$.
There is a quotient $\gg_1^\prime$ of $\gg_1$ such that $\hg(Y) = \gg_1^\prime
\oplus \gg_3$. Notice that via $\hg(Y) \leftisomarrow \hg(Y \times X_2) =
\gg_1^\prime \oplus \gg_3 \twoheadrightarrow \hg(X_2)$ we obtain an isomorphism
$\hg(X_2) \isomarrow \gg_3$ mapping $J_{X_2}$ to the $\gg_3$-component of
$J_Y$. Lemma (\refLemIdeal) then gives $\Hom(X_2,Y) \neq 0$. $\square$
\whiteB

\noindent
\advance\secno by 1
{\bf (\chapno.\sectno) Remark.} It was shown by Borovoi [\refBorov] that
$\hg(X)$ is $\mQ$-simple if $\Eo(X) = \mQ$. For a generalization of this result
to absolutely irreducible Hodge structures of arbitrary level see [\refZarLNM].
\whiteB

\noindent
\advance\secno by 1
\edef\refLemProdTor{\number\chno.\number\secno}%
{\bf (\chapno.\sectno) Lemma.}\enspace {\sl Let $X_1$ and $X_2$ be nonzero
complex abelian varieties. Assume that the Hodge group $\Hg(X_2)$ is a
$\mQ$-simple algebraic torus. (In particular $X_2$ is of CM-type.) Write $X =
X_1 \times X_2$. If $\Hg(X) \neq \Hg(X_1) \times
\Hg(X_2)$ then the center of $\Hg(X_1)$ contains an algebraic torus which is
$\mQ$-isogenous to $\Hg(X_2)$.}
\whiteD

\Proof~Suppose that $\Hg(X) \neq \Hg(X_1) \times \Hg(X_2)$. The assumption that
$\Hg(X_2)$ is $\mQ$-simple implies that $\hg(X_2)$ does not contain a proper
algebraic Lie subalgebra. Using the notations of (\refProd) we then have that
$\hg(X) = \gg_1 \oplus \gg_3 \isomarrow \hg(X_1)$ and $\hg(X_2) \cong \gg_3$.
This readily implies the lemma, noting that $\gg_1$ and $\gg_3$ are algebraic
Lie subalgebras of $\hg(X)$. $\square$
\whiteC

Next let us recall a lemma from [\refLZTate] that was also used in
[\refMZDuke]. This lemma was also stated in [\refHazCM], where it is attributed
to Ribet. To formulate it, we need the following notation. Suppose $F$ is a
CM-field containing an imaginary quadratic field $k$. In \S 2 above we defined
the algebraic torus $\UU_F$ over $\mQ$. The subfield $k \subset F$ gives rise
to a subtorus $\SU_{F/k} \subset \UU_F$ of codimension 1, by
$$
\SU_{F/k} = \Ker(\Nm_{F/k} \colon \UU_F \rightarrow \UU_k)\, .
$$
With this notation, we have the following lemma. For a proof we refer to
[\refMZDuke].
\whiteB

\noindent
\advance\secno by 1
\edef\refDukeLem{\number\chno.\number\secno}%
{\bf (\chapno.\sectno) Lemma.}\enspace {\sl Let $F$ be a CM-field. Suppose $H$
is an algebraic subtorus of $\UU_F$ of codimension 1. Then there exists an
imaginary quadratic subfield $k \subset F$ such that $H = \SU_{F/k}$.}
\whiteC

Combining the above lemmas with the facts in (\refEllC) gives the following
result.
\whiteB

\noindent
\advance\secno by 1
\edef\refEllFact{\number\chno.\number\secno}%
{\bf (\chapno.\sectno) Proposition.}\enspace {\sl Let $X$ be an abelian variety
and let $E$ be an elliptic curve, both over $\mC$. Suppose $\Hom(E,X) = 0$.
Then either $\Hg(X \times E) = \Hg(X) \times \Hg(E)$ or $\Eo(E) = k$ is an
imaginary quadratic field such that there exists an embedding of $k$ into the
center of $\Eo(X)$.}
\whiteD

\Proof~If $\Eo(E) = \mQ$ then we apply Lemma (\refLemProdSS). Hence we may
assume that $\Eo(E) = k$ is an imaginary quadratic field, so that $\Hg(E)$ is
the rank 1 torus $\UU_k$.

Write $C$ for the center of $\Eo(X)$. Then $C$ has the form $C = K_1 \times
\cdots \times K_m \times F_1 \times \cdots \times F_n$, where $K_1, \ldots,
K_m$ are totally real fields and $F_1, \cdots, F_n$ are CM-fields. The center
$Z$ of $\Hg(X)$ is contained in $\UU_{F_1} \times \cdots \times \UU_{F_n}$. By
Lemma (\refLemProdTor), if $\Hg(X \times E) \neq \Hg(X) \times \Hg(E)$ then
there is a homomorphism $\UU_k \rightarrow \UU_{F_1} \times \cdots \times
\UU_{F_n}$ with finite kernel. If $\UU_{F_i}$ is a factor such that the
projection of $\UU_k$ to $\UU_{F_i}$ has rank 1 then it easily follows from
Lemma (\refDukeLem) that there exists an embedding $k \rightarrow F_i$. This
proves the claim. $\square$
\whiteC

As an easy corollary we obtain a result first proven by Imai [\refImai].
\whiteB

\noindent
\advance\secno by 1
\edef\refImaiThm{\number\chno.\number\secno}%
{\bf (\chapno.\sectno) Corollary.}\enspace {\sl Let $X_1, \ldots, X_n$ be
elliptic curves over $\mC$, no two of which are isogenous. Write $X = X_1
\times \cdots \times X_n$. Then $\Hg(X) = \Hg(X_1) \times \cdots \times
\Hg(X_n)$. In particular, every product of elliptic curves satisfies condition
{\rm (D)} in {\rm (\refHazMur)}.}
\whiteD

\Proof~Immediate from the proposition, by induction on the number of factors.
$\square$
\whiteB

\noindent
\advance\secno by 1
{\bf (\chapno.\sectno) Remark.} The constructions in this section were inspired
by similar results for abelian varieties over finite fields obtained in
[\refZarTaFin].
\whiteB\whiteB

\advance\chno by 1
\secno=0
\noindent
\centerline{{\bf \S \chapno.~Hodge groups of simple abelian surfaces of
CM-type.}}
\whiteA

\noindent
In this section we study Hodge groups of simple abelian surfaces of CM-type. We
use this to prove Theorem~(\refMainThm) for the product of two such surfaces.
\whiteB

\noindent
\advance\secno by 1
\edef\refCM{\number\chno.\number\secno}%
{\bf (\chapno.\sectno)}\enspace Let $F$ be a CM-field. Write $\Sigma_F$ for the
set of embeddings $F \rightarrow \mC$. Let $\iota \colon x \mapsto \bar{x}$
denote the complex conjugation on $F$. (Recall that $\iota$ is independent of
the choice of an embedding of $F$ into $\mC$.) By a CM-type for $F$ we mean a
subset $\Phi \subset \Sigma_F$ such that, writing $\overline{\Phi} =
\{\bar{\phi}\mid \phi \in \Phi\}$, we have $\Sigma_F = \Phi \amalg
\overline{\Phi}$.

Write $F_0 \subset F$ for the totally real subfield. The choice of a CM-type
$\Phi$ for $F$ is equivalent to giving an identification $F \otimes_\mQ \mR
\isomarrow \mC^{\Sigma_{F_0}}$. Writing $J = J_\Phi \in F \otimes_\mQ \mR$ for
the element which maps to $(i,i, \ldots,i)$ we obtain a bijection
$$
\{\hbox{CM-types for $F$}\} \isomarrow \Gamma_F := \{J \in F \otimes_\mQ \mR
\mid J^2 = -1\}
$$
which is equivariant for the natural $\Aut(F)$-action on both sides.

To the CM-type $(F,\Phi)$ we can associate an isogeny class of complex abelian
varieties by taking $F$ as a $\mQ$-lattice and $J_\Phi$ as a complex structure.
Two CM-types $(F,\Phi)$ and $(F,\Psi)$ give rise to the same isogeny class if
and only if there exists an automorphism $\alpha \in \Aut(F)$ with $\Psi =
{}^\alpha\Phi$. Note that if $X$ is an abelian variety in the isogeny class
associated to $(F,\Phi)$ then $J_\Phi$ is just the operator $J_X$ as in
(\refHgLie). We have $J_{\overline{\Phi}} = - J_\Phi$.

Now let $F$ be a quartic CM-field which does not contain an imaginary quadratic
subfield. Then either (\romannumeral1)~$F$ is Galois over $\mQ$, in which case
$\Aut(F)$ is cyclic of order 4 acting transitively on $\Gamma_F$, or
(\romannumeral2)~$F$ is not Galois over $\mQ$, its normal closure $L$ has
degree 8 over $\mQ$, and $\Aut(F) = \{\id,\iota\}$. In case (\romannumeral1)
there is only one isogeny class of abelian surfaces with CM by $F$, in case
(\romannumeral2) there are two such isogeny classes.
\whiteB

\noindent
\advance\secno by 1
\edef\refCMAbSurfCor{\number\chno.\number\secno}%
{\bf (\chapno.\sectno) Proposition.}\enspace {\sl Let $X_1$ and $X_2$ be two
simple abelian surfaces with CM by the same quartic CM-field $F$. Suppose $X_1$
and $X_2$ are not isogenous. Write $X = X_1 \times X_2$. Then $\Hg(X) =
\Hg(X_1) \times \Hg(X_2)$.}
\whiteD

\Proof~Fix isomorphisms $F \cong \Eo(X_i)$; this gives identifications
$\Hg(X_i) = \UU_F$. As just explained, the assumption that $X_1 \not\sim X_2$
implies that $F$ is not Galois over $\mQ$. {\sl A priori\/} the Galois group
$\Gal(L/\mQ)$ could be isomorphic to either the dihedral group $D_4$ or the
quaternion group $Q$. By [\refShimFam], Propositions 14 and 18, the CM-field
$F$ does not contain an imaginary quadratic field. Lemma (\refDukeLem) then
shows that the torus $\UU_F$ is $\mQ$-simple. Its splitting field is the field
$L$, as one verifies without great difficulty. Writing $X^\ast = X^\ast(\UU_F)$
for the character group, the previous facts mean that $X^\ast_\mQ$ is a
faithful irreducible
2-dimensional $\mQ$-representation of $\Gal(L/\mQ)$. Now remark that the group
$Q$ does not
admit such a representation (cf.\ [\refSerRGF], Sect.\ 12.2, p.\ 108). Hence
$\Gal(L/\mQ) \cong D_4$.

Consider the ``standard'' representation $\rho
\colon D_4 \rightarrow \GL_{2}(\mQ)$, realizing $D_4$ as the subgroup of
$\GL_2(\mZ)$ generated by the matrices
$$
\left(\matrix{0&-1\cr1&0\cr}\right)
\qquad \hbox{and}\qquad
\left(\matrix{0&1\cr1&0}\right)\, .
$$
We remark that this $\rho$ is the {\sl only\/} faithful irreducible
2-dimensional $\mQ$-representation of $D_4$ (up to isomorphism), and that it is
absolutely irreducible. (This is an elementary exercise.) This last fact
implies that $\End(\UU_F)=\mZ$, i.e., all homomorphisms $\UU_F \to \UU_F$ are
of the form $x \mapsto x^m$ for some integer $m$.

Assume that $\Hg(X) \neq \Hg(X_1) \times \Hg(X_2)$. The fact that $\UU_F$ is
$\mQ$-simple  implies that both surjective projection maps $\Hg(X) \to
\Hg(X_i)$ are isogenies and therefore $\Hg(X)$ is isogenous to $\UU_F$. This
implies that
 $\Hom(\Hg(X),\UU_F)$ is isomorphic to $\mZ$ as an abelian group.
Let $u\colon \Hg(X) \to \UU_F$ be a generator of this group. Clearly $u$ is an
isogeny. The homomorphism $j \colon \Hg(X) \rightarrow \Hg(X_1) \times \Hg(X_2)
= \UU_F \times \UU_F$ is of the form $x \mapsto (u(x)^m,u(x)^n)$ for some
integers $m$ and $n$. In particular, $\Ker(u) \subseteq \Ker(j)$. As $j$ is
injective, it follows that $u$ is an isomorphism and that the integers $m$ and
$n$ are relatively prime.

Under $\pr_i$, the element $J_X \in \hg(X) \otimes \mC$ is mapped to $J_i =
J_{X_i}$. Under the given identifications $\hg(X_1) = \uu_F = \hg(X_2)$ we thus
find that $J_2 = (n/m) \cdot J_1$. Since both $J_1$ and $J_2$, viewed as
elements of $F$, satisfy $J_i^2 = -1$ it follows that $m = \pm n$. But $m$ and
$n$ are relatively prime, so $m,n \in \{\pm 1\}$ and $J_1 = \pm J_2$. This
implies that $X_1$ and $X_2$ are isogenous (see (\refCM)), contradicting the
assumptions. $\square$
\whiteB\whiteB

\advance\chno by 1
\secno=0
\noindent
\centerline{{\bf \S \chapno.~Proof of the main result.}}
\whiteA

\noindent
\advance\secno by 1
{\bf (\chapno.\sectno)}\enspace Let $X$ be a complex abelian variety with $g =
\dim(X) \leq 4$. Our first goal is to prove (\refMainThm). As recalled above we
already know this in case $X$ is simple. In the rest of this section we may,
and will, therefore assume that $X$ is {\sl not simple}.

Up to isogeny we can decompose $X$ as $X \sim Y_1^{m_1} \times \cdots \times
Y_r^{m_r}$ where $Y_1, \ldots,Y_r$ ($r \in \mZ_{\geq 1}$) are simple, pairwise
non-isogenous abelian varieties and $m_1, \ldots, m_r \in \mZ_{\geq 1}$.
Correspondingly, the endomorphism algebra $D$ decomposes as $D = D_1 \times
\cdots \times D_r$ where $D_i = \Eo(Y_i^{m_i}) \cong M_{m_i}(\Eo(Y_i))$. Write
$V = H_1(X,\mQ)$ and $V_i = H_1(Y_i^{m_i},\mQ)$. Choose polarizations
$\lambda_i$ of $Y_i^{m_i}$, let $\lambda$ be the ``product'' polarization
$\lambda = (\lambda_1,\ldots,\lambda_r)$ of $X$, and let $\phi_i \colon V_i
\times V_i \rightarrow \mQ$ resp.\ $\phi \colon V \times V \rightarrow \mQ$ be
the associated Riemann forms. With these notations we have the obvious remark
that $\Hg(X) = \Sp_D(V,\phi)$ if and only if $\Hg(X) = \Hg(Y_1^{m_1}) \times
\cdots \times \Hg(Y_r^{m_r})$ and $\Hg(Y_i^{m_i}) = \Sp_{D_i}(V_i,\phi_i)$ for
all $i$.

Now assume that $X$ is not simple with $\dim(X) = 4$. Note that $X$ has no
factors of Type \Romno 3 (since Type \Romno 3 does not occur in dimension $\leq
3$). Case (a) of the introduction will be dealt with in (5.3) below. If we are
not in case (a) then, using the Theorem (\refHazMur) of Hazama and Murty and
the results discussed in \S 2, we see
that in order to prove (\refMainThm) for $X$ it suffices to show that $\Hg(X) =
\Hg(Y_1^{m_1}) \times \cdots \times \Hg(Y_r^{m_r})$.
\whiteB

\noindent
\advance\secno by 1
{\bf (\chapno.\sectno)}\enspace Suppose $g =3$. Suppose also that $X$
decomposes, up to isogeny, as a product $X \sim X_1 \times X_2$ of an elliptic
curve $X_1$ and a simple abelian surface $X_2$. Then the center of $\Eo(X_2)$
does not contain an imaginary quadratic field. By Proposition~(\refEllFact) it
follows that $\Hg(X) = \Hg(X_1) \times \Hg(X_2)$.

Combining this with Corollary~(\refImaiThm), we have proven (\refMainThm) in
case
$\dim(X) \leq 3$. In particular, for every complex abelian variety $X$ of
dimension $\leq 3$ we have $\Hg(X) = \Sp_D(V,\phi)$ and condition (D) in
(\refHazMur) is satisfied.
\whiteB

\noindent
\advance\secno by 1
\edef\refArgCasea{\number\chno.\number\secno}%
{\bf (\chapno.\sectno)}\enspace Let $X$ be a complex abelian variety which is
isogenous to a product, say $X \sim X_1 \times X_2$, where $X_1$ is an elliptic
curve and $X_2$ is a simple abelian threefold. Suppose furthermore that $k :=
\Eo(X_1)$ is an imaginary quadratic field and that there exists an embedding $k
\hookrightarrow F := \Eo(X_2)$. This means we are in case (a) of the
introduction. Either (a1) $F = k$, or (a2) $F$ is a sextic CM-field.

Embed $k$ as a subfield of $\Eo(X)$ such that it acts with multiplicities
$(2,2)$ on the tangent space $T_{X,0}$. (Our assumption that $X_2$ is simple
implies that $k$ acts on $T_{X_2,0}$ with multiplicities $(1,2)$, see
[\refShimFam], Proposition 14. Therefore, if we fix $\Eo(X_1) = k
\hookrightarrow
\Eo(X_2)$ then either $\alpha \mapsto (\alpha, \alpha) \in \Eo(X_1) \times
\Eo(X_2)$ or $\alpha \mapsto (\bar{\alpha},\alpha)$ gives an embedding as
required.) Then the space $W_k \subset H^4(X,\mQ)$ consists of Hodge classes.
We know that
$$
\Hg(X) \subseteq \Hg(X_1) \times \Hg(X_2) = \cases{\UU_k
\times \UU_k(V_{X_2},\psi_{X_2})&in case (a1);\cr\noalign{\smallskip}\UU_k
\times \UU_F&in case (a2).\cr}
$$
(See \S 2 for notations.) The Hodge group
acts trivially on $W_k$, i.e., its elements have trivial $k$-linear
determinant. We then easily find that we must have
$$
\Hg(X) = \cases{\{(u_1,u_2) \in \UU_k \times \UU_k(V_{X_2},\psi_{X_2}) \bigm|
u_1
\cdot \det_k(u_2) = 1\}&in case (a1);\cr\noalign{\smallskip}\{(u_1,u_2) \in
\UU_k \times \UU_F \bigm| u_1 \cdot \det_k(u_2) =
1\}&in case (a2),}
$$
where $\det_k \colon \UU_k(V_{X_2},\psi_{X_2}) \rightarrow \UU_k$ denotes the
$k$-linear determinant map, resp.\ $\det_k = \Nm_{F/k} \colon \UU_F \rightarrow
\UU_k$. (To see our claim, note that $\UU_k$ has rank 1 and that $\Hg(X)$ maps
surjectively onto $\Hg(X_2)$, so $\Hg(X)$ can at most have codimension 1 in
$\Hg(X_1) \times \Hg(X_2)$.)

The K\"unneth decomposition gives
$$
\eqalign{H^4(X,\mQ) = &[H^2(X_1,\mQ) \otimes H^2(X_2,\mQ)]\;\oplus\;
[H^1(X_1,\mQ)
\otimes H^3(X_2,\mQ)]\cr&\oplus\; [H^0(X_1,\mQ) \otimes H^4(X_2,\mQ)]\, .}
$$
The Hodge classes in $H^2(X_1,\mQ) \otimes H^2(X_2,\mQ) \cong H^2(X_2,\mQ)(-1)$
and those in $H^0(X_1,\mQ) \otimes H^4(X_2,\mQ) \cong H^4(X_2,\mQ)$ are linear
combination of products of divisor classes. The space of Weil classes $W_k$ is
a subspace of $H^1(X_1,\mQ) \otimes H^3(X_2,\mQ)$. (Since we are viewing $W_k$
as a {\sl subspace\/} of $H^4(X,\mQ)$, rather than a {\sl quotient\/}, some of
our identifications may seem a little unnatural, cf.\ [\refMZCrel], Sect.~7.)

We have an isomorphism of Hodge structures
$$
H^1(X_1,\mQ) \otimes H^3(X_2,\mQ) \cong \Hom(H^1(X_1,\mQ),H^3(X_2,\mQ))(-1)\, .
$$
Under this isomorphism, the space of Hodge classes in $H^1(X_1,\mQ) \otimes
H^3(X_2,\mQ)$ corresponds to the space
$\Hom_\HS(H^1(X_1,\mQ),H^3(X_2,\mQ))(-1)$ of homomorphisms of $\mQ$-Hodge
structures. The Hodge structure $H^1(X_1,\mQ)$ is irreducible and has
endomorphism ring $k$. Therefore, our assertion that $W_k$ is the space of
Hodge classes in $H^1(X_1,\mQ) \otimes H^3(X_2,\mQ)$ is equivalent to saying
that $H^3(X_2,\mQ)$ contains only one
copy of $H^1(X_1,\mQ)(-1)$ as a rational sub-Hodge structure. It suffices to
prove this in case (a2), since the group $\UU_k(V_{X_2},\psi_{X_2})$ contains
tori of the form $\UU_F$
where $F$ is a sextic CM-field containing $k$. (Put differently: we can
specialize from case (a1) to case (a2).)

Suppose then that $W \subset H^3(X_2,\mQ)$ is a rational sub-Hodge structure
isomorphic to
$H^1(X_1,\mQ)(-1)$. Then $\Hg(X_2) = \UU_F$ acts on $W$ through the torus
$\Hg(W) = \Hg(X_1) = \UU_k$. The kernel of the corresponding homomorphism
$\UU_F \rightarrow \UU_k$ is necessarily the subtorus $\SU_{F/k} \subset
\UU_F$. (Cf.\ Lemma (\refDukeLem).) But now we remark that the space of
$\SU_{F/k}$-invariants in $H^3(X_2,\mQ)$ has $\mQ$-dimension 2, which proves
our claim. (In fact, the space of $\SU_{F/k}$-invariants in $H^3(X_2,\mQ)$ is
precisely the subspace $W_k(X_2) \subset H^3(X_2,\mQ)$, which is naturally a
1-dimensional $k$-vector space.)

In sum, the previous arguments prove (\refMainThm) for case (a).
\whiteB

\noindent
\advance\secno by 1
{\bf (\chapno.\sectno)}\enspace Let $X$ be a non-simple complex abelian
fourfold. Suppose $X$ is not of CM-type. Then $X$ contains a simple abelian
subvariety $X_2$ which is not of CM-type. We can write $X \sim X_1 \times
X_2^r$ with $r \geq 1$ and $\Hom(X_2,X_1) = 0$.

Suppose that we are not in case (a). We want to show that $\Hg(X) = \Hg(X_1)
\times \Hg(X_2^r)$. If $r > 1$ then we are reduced to the case $g \leq 3$,
since $\Hg(X_1 \times X_2^r) \cong \Hg(X_1 \times X_2)$. Assume then that $r
=1$. We distinguish two cases. If $\dim(X_2) = 3$ then $X_1$ is an elliptic
curve and we can apply (\refEllFact), which works since we are not in case (a).
If $\dim(X_2) < 3$ then the desired equality $\Hg(X) = \Hg(X_1) \times
\Hg(X_2)$ follows from (\refPropLowg) and (\refLemProdSS).
\whiteB

\noindent
\advance\secno by 1
{\bf (\chapno.\sectno)}\enspace Let $X$ be a non-simple complex abelian
fourfold of CM-type. Suppose that $X$ is isogenous to $X_1 \times X_2^r$ with
$\dim(X_1) = 1$ and $\Hom(X_1,X_2) = 0$. If we are not in case (a) then there
is no embedding of $\Eo(X_1)$ into the center of $\Eo(X_2)$. It thus follows
from Proposition~(\refEllFact) that $\Hg(X) = \Hg(X_1) \times \Hg(X_2^r)$.

This only leaves us with the case where $X \sim X_1 \times X_2$, with $X_1$ and
$X_2$ simple abelian surfaces. If $X_1$ and $X_2$ are isogenous then we are
done. If $X_1$ and $X_2$ are not isogenous then Proposition (\refCMAbSurfCor)
shows that $\Hg(X) = \Hg(X_1) \times \Hg(X_2)$.
\whiteC

This completes the proof of Theorem (\refMainThm).
\whiteB

We now turn to the proof of Theorem~(\refThmTwo). As we have seen in \S 2, the
statement is known if $X$ is simple. So again we may, and shall, assume $X$ to
be non-simple. Furthermore we can assume that every simple factor of $X$ occurs
with multiplicity 1.

Write $X_1 \subseteq X$ for the maximal abelian subvariety which has no factors
of Type \Romno 4, and $X_2 \subseteq X$ for the maximal abelian subvariety of
which all factors are of Type \Romno 4. Write $d_i = \dim(X_i)$. We shall treat
the possibilities case by case.
\whiteB

\noindent
\advance\secno by 1
{\bf (\chapno.\sectno)}\enspace Suppose $(d_1,d_2) = (5,0)$, so that $X$ has no
factors of Type \Romno 4. If $X$ contains an elliptic curve $E$ then $\Eo(E) =
\mQ$ (since $E$ is not of Type \Romno 4) and Theorem~(\refThmTwo) follows by
Proposition~(\refEllFact). If $X$ does not contain an elliptic curve then all
its
simple factors satisfy condition (D) in (\refHazMur) and we conclude using
Theorem~(\refHazThm).
\whiteB

\noindent
\advance\secno by 1
{\bf (\chapno.\sectno)}\enspace Suppose $(d_1,d_2) = (4,1)$ or $(d_1,d_2) =
(3,2)$. Then $X_2$ is of CM-type and Theorem~(\refHazThm) gives $\Hg(X) =
\Hg(X_1) \times \Hg(X_2)$.
\whiteB

\noindent
\advance\secno by 1
{\bf (\chapno.\sectno)}\enspace Suppose $(d_1,d_2) = (2,3)$. If $X_1$ is simple
then Lemma (\refLemProdSS) gives $\Hg(X) = \Hg(X_1) \times \Hg(X_2)$. If $X_1$
is not simple then it is isogenous to a product of two elliptic curves, $X_1
\sim E_1 \times E_2$ with $\Eo(E_1) = \Eo(E_2) = \mQ$ and where we may assume
$E_1$ and $E_2$ to be non-isogenous. Proposition~(\refEllFact) then gives
$\Hg(X) =
\Hg(E_1) \times \Hg(E_2) \times \Hg(X_2)$.
\whiteB

\noindent
\advance\secno by 1
\edef\refOneFour{\number\chno.\number\secno}%
{\bf (\chapno.\sectno)}\enspace Suppose $(d_1,d_2) = (1,4)$. Then $\Hg(X) =
\Hg(X_1) \times \Hg(X_2)$ by Proposition~(\refEllFact).
\whiteB

\noindent
\advance\secno by 1
\edef\refNilFive{\number\chno.\number\secno}%
{\bf (\chapno.\sectno)}\enspace From now on, let us assume that $(d_1,d_2) =
(0,5)$, meaning that all simple factors of $X$ are of Type \Romno 4. Let
$d_{\min}$ be the minimal dimension of a simple factor of $X$. Since we assume
$X$ to be non-simple we have $d_{\min} = 1$ or $d_{\min} = 2$.

First suppose that $d_{\min} = 2$. Then $X \sim Y_1 \times Y_2$ where $Y_1$ is
a simple abelian surface and $Y_2$ is a simple abelian threefold. Note that
$Y_1$ is of CM-type with $\Hg(Y_1) = \UU_{F_1}$, where $F_1 = \Eo(Y_1)$.

If $Y_2$ is not of CM-type then Lemma (\refLemProdTor) readily gives $\Hg(X) =
\Hg(Y_1) \times \Hg(Y_2)$. If $Y_2$ is of CM-type then $F_2 = \Eo(Y_2)$ is a
sextic CM-field. By Lemma (\refLemProdTor) and Lemma (\refDukeLem), we can have
$\Hg(X) \neq \Hg(Y_1) \times \Hg(Y_2)$ only if $F_2$ contains an imaginary
quadratic field $k$ such that $\UU_{F_1}$ is isogenous to $\SU_{F_2/k}$.
Suppose this is the case. Write $\Omega_1$ for the normal closure of $F_1$ over
$\mQ$. Either $\Omega_1 = F_1$ and $\Gal(\Omega_1/\mQ) = \mZ/4\mZ$ (as $F_1$
does not contain an imaginary quadratic field), or $\Omega_1$ has degree 8 over
$\mQ$. Next write $K_2$ for the totally real subfield of $F_2$ and let
$\Omega_2$ be the normal closure of $K_2$ over $\mQ$. As $F_2$ contains the
imaginary quadratic field $k$, the normal closure of $F_2$ over $\mQ$ is the
compositum $k \cdot \Omega_2$. The Galois group $\Gal(k \cdot \Omega_2/\mQ)$ is
either $(\mZ/2\mZ) \times {\frak S}_3$ or $(\mZ/2\mZ) \times (\mZ/3\mZ)$. Now
$\Omega_1$ is the splitting field of the $\mQ$-torus $\UU_{F_1}$ and $k \cdot
\Omega_2$ contains the splitting field of $\UU_{F_2}$. The assumption that
$\UU_{F_1}$ is isogenous to $\SU_{F_2/k}$ thus implies that $\Omega_1 \subseteq
k \cdot \Omega_2$. Looking at Galois groups we obtain a contradiction. Hence
again $\Hg(X) = \Hg(Y_1) \times \Hg(Y_2)$.
\whiteB

\noindent
\advance\secno by 1
{\bf (\chapno.\sectno)}\enspace From now on, let us assume that $(d_1,d_2) =
(0,5)$ and that $d_{\min} = 1$. Write $X \sim E \times Y$, where $E$ is an
elliptic curve and $\dim(Y) = 4$. Without loss of generality we may assume that
$\Hom(E,Y) = 0$. (If not then we are reduced to the case $\dim(X) \leq 4$.) Let
$d_{\max}$ be the maximal dimension of a simple factor of $X$.

If $d_{\max} \leq 2$ then all simple factors of $Y$ are of CM-type and there
does not exist an embedding of $\Eo(E)$ into the center of $\Eo(Y)$. Then
Proposition~(\refEllFact) gives $\Hg(X) = \Hg(E) \times \Hg(Y)$.

If $d_{\max} = 3$ then $Y$ is isogenous to a product of an elliptic curve $Y_1$
and a simple abelian threefold $Y_2$. If $\Eo(Y_2)$ contains an imaginary
quadratic field then this subfield is unique. Therefore, possibly after
interchanging the roles of $E$ and $Y_1$ we find that there does not exist an
embedding of $\Eo(E)$ into the center of $\Eo(Y)$. (Note that $\Eo(E) =
\Eo(Y_1)$ implies that $E \sim Y_1$, which we excluded.) Again by
Proposition~(\refEllFact) we then find $\Hg(X) = \Hg(E) \times \Hg(Y)$.

Finally, let us assume that $d_{\max} = 4$, i.e., that $Y$ is simple (of Type
\Romno 4). Write $k = \Eo(E)$ and $F = \Eo(Y)$. If there is no embedding $j
\colon k \hookrightarrow F$ then Proposition~(\refEllFact) gives $\Hg(X) =
\Hg(E)
\times \Hg(Y)$. Suppose then that there exists an embedding $j$. We distinguish
2 cases.

{\sl Case 1:\/} Suppose that $k$ acts on $T_{Y,0}$ with multiplicities $(2,2)$,
so that $\Hg(Y) = \SU_{F/k}$. Suppose also that $\Hg(X) \neq \Hg(E) \times
\Hg(Y)$. As $\Hg(E)$ has rank 1 we find that the (surjective) homomorphism
$\pr_2 \colon \Hg(X) \rightarrow \Hg(Y)$ is an isogeny. Then there also exists
an isogeny $f \colon \SU_{F/k} = \Hg(Y) \rightarrow \Hg(X)$ and we obtain a
non-trivial homomorphism $\alpha := \pr_1 \circ f \colon \SU_{F/k} \rightarrow
\UU_k = \Hg(E)$. Next choose a homomorphism $j \colon \UU_F \rightarrow
\SU_{F/k}$ such that the composition $\SU_{F/k} \hookrightarrow \UU_F
\mapright{j} \SU_{F/k}$ is an isogeny. The identity
component $K$ of $\Ker(\alpha \circ j \colon \UU_F \rightarrow \UU_k)$ is a
codimension 1 subtorus of $\UU_F$. By (\refDukeLem) there exists an imaginary
quadratic subfield $l \subset F$ such that $K = \SU_{F/l}$. The quotient
$U_F/K$ is isogenous to $\UU_l$ but also to $\UU_k$ (as $K = \Ker(\alpha \circ
j)$). Hence $\UU_k$ is isogenous to $\UU_l$, which implies that $k = l$ and $K
= \SU_{F/k}$. (Note that $F$ has only one subfield isomorphic to $k$.) It is
clear though from our construction that $\SU_{F/k}$ is not contained in
$\Ker(\alpha \circ j)$, as $\pr_1$ is surjective. Hence $\Hg(X) = \Hg(E) \times
\Hg(Y)$.

{\sl Case 2:\/} Suppose that $k$ acts on $T_{Y,0}$ with multiplicities $(1,3)$.
(By [\refShimFam], Proposition 14 this is the only other case that occurs.)
Rather
than looking at $E \times Y$, let us look at $Z := E^2 \times Y$. There is an
embedding $k \hookrightarrow \Eo(Z)$ such that $k$ acts on $T_{Z,0}$ with
multiplicities $(3,3)$. This implies that the corresponding space of Weil
classes $W_k \subset H^6(Z,\mQ)$ consists of Hodge classes and that $\Hg(Z)
\subseteq \SU_k(V_Z,\psi)$. (For this last conclusion, see [\refMZDuke], Lemma
2.8.) Returning to our original abelian variety $X \sim E \times Y$ we find
that $\Hg(X)$ is contained in the subgroup $H \subset \Hg(E) \times \Hg(Y) =
\UU_k \times \UU_F(V_Y,\psi)$ given by
$$
H = \{(u_1,u_2) \in \UU_k \times \UU_F(V_Y,\psi) \mid u_1^2 \cdot \det_k(u_2) =
1\}\, ,
$$
where $\det_k \colon \Hg(Y) = \UU_F(V_Y,\psi) \rightarrow \UU_k$ is the
$k$-linear determinant. Now remark that $\UU_k$ has rank 1, so that the
projection $H \rightarrow \Hg(Y)$ is an isogeny. As $H$ is connected and $\pr_2
\colon \Hg(X) \rightarrow \Hg(Y)$ is surjective, we conclude that $\Hg(X) = H$.
\whiteB

\noindent
\advance\secno by 1
{\bf (\chapno.\sectno)}\enspace We have now computed the Hodge groups of all
complex abelian 5-folds. It remains to be shown that this indeed gives the
conclusions as stated in Theorem~(\refThmTwo). Part (\romannumeral4) of the
theorem follows by going through the above and using (\refMainThm) and
(\refPrimeDim). All that remains to be done is the computation of the Hodge
rings in the cases (e), (f) and (g).

Case (f) is easy. It was established in (\refOneFour) and (\refNilFive) that
$\Hg(X) = \Hg(X_0) \times \Hg(X_1 \times X_2)$. (Notations as in the
introduction.) The rest of statement (\romannumeral2) of (\refThmTwo) readily
follows.

Next suppose we are in case (e). By the duality $H^j(X,\mQ)(5) \cong
H^{10-j}(X,\mQ)^\vee$ we only have to show that $\cB^2(X) \subset H^4(X,\mQ)$
is generated by $\cD^2(X)$ and the spaces $W_{k,\alpha}$. The K\"unneth formula
gives
$$
\eqalign{H^4(X,\mQ) = &\phantom{\oplus\;\ } [H^4(X_1^2,\mQ) \otimes
H^0(X_2,\mQ)]\;\oplus\;
[H^3(X_1^2,\mQ) \otimes H^1(X_2,\mQ)]\cr
&\oplus\; [H^2(X_1^2,\mQ) \otimes H^2(X_2,\mQ)]\;\oplus\; [H^1(X_1^2,\mQ)
\otimes H^3(X_2,\mQ)]\cr
&\oplus\; [H^0(X_1^2,\mQ) \otimes H^4(X_2,\mQ)]\, .\cr
}
$$
In $H^4 \otimes H^0$ and $H^0 \otimes H^4$ we only have decomposable classes.
In
$$
H^3(X_1^2,\mQ) \otimes H^1(X_2,\mQ) \cong \Hom(H^1(X_1^2,\mQ),
H^1(X_2,\mQ))(-2)
$$
there are no non-zero Hodge classes, as there are no non-zero homomorphisms
from $X_1^2$ to $X_2$. Next we have $H^1(X_1^2,\mQ) \otimes H^3(X_2,\mQ) \cong
[H^1(X_1,\mQ) \otimes H^3(X_2,\mQ)]^{\oplus 2}$, so that the Hodge classes in
$H^1 \otimes H^3$ are just the elements of the spaces $W_{k,\alpha}$. (See
(\refArgCasea). Also note that in fact we only need two spaces $W_{k,\alpha_1}$
and $W_{k,\alpha_2}$ for ``linear independent'' choices $\alpha_1$ and
$\alpha_2$.)

To settle case (e) it thus remains to compute the Hodge classes in $H^2 \otimes
H^2$. Write $V_2 = H_1(X_2,\mQ)$ and $F = \Eo(X_2)$. Either $F = k$ or $F$ is a
sextic CM-field. Fix an element $a \in F$ with $\bar{a} = -a$. The Hodge group
$\Hg(X_2)$ is the unitary group $\UU_F(V_2,\psi)$, where $\psi \colon V_2
\times
V_2 \rightarrow F$ is an $F$-hermitian form such that $\trace_{F/\mQ}(a \cdot
\psi)$ is the Riemann form of a polarization. (See (\refThreef) and notice that
if $F$ is a sextic CM-field then the same description works, since
$\UU_F(V_2,\psi)$ in that case is just the torus $\UU_F$.)

Consider the algebraic $\mQ$-subgroup $\SU_{F/k}(V_2,\psi) = \Ker(\det_k \colon
\UU_F(V_2,\psi) \rightarrow \UU_k)$. We claim that $\SU_{F/k}(V_2,\psi)$ and
$\UU_F(V_2,\psi)$  have the same centralizer in $\End(V_2)$. To see this we can
extend scalars from $\mQ$ to $\mC$ and consider the actions of
$\SU_{F/k}(V_2,\psi) \otimes_\mQ \mC$ and $\UU_F(V_2,\psi) \otimes_\mQ \mC$ on
$V_2 \otimes_\mQ \mC$. Treating the cases $F=k$ and $[F:\mQ] = 6$ separately,
the claim is then easily verified. As $H^2(X_2,\mQ)$ is isomorphic to a
sub-Hodge structure
of $\End(V_2)(-1)$ it follows that the space of
$\SU_{F/k}(V_2,\psi)$-invariants
in $H^2(X_2,\mQ)$ is equal to the space $\cB^1(X_2)$ of $\Hg(X_2)$-invariants.
Now our description of $\Hg(X) \cong \Hg(X_1 \times X_2)$ in (\refArgCasea)
above shows that $\Hg(X) \supset \{1\} \times \SU_{F/k}(V_2,\psi)$, so that the
Hodge classes in $H^2 \otimes H^2$ are contained in $H^2(X_1^2,\mQ) \otimes
\cB^1(X_2)$. But then it readily follows that the Hodge classes in $H^2 \otimes
H^2$ all lie in $\cB^1(X_1^2) \otimes \cB^1(X_2)$ and are therefore
decomposable. This finishes the proof of (\romannumeral1) of
Theorem~(\refThmTwo).

Finally, suppose we are in case (g). Again we only have to look at
$H^4(X,\mQ)$. The only interesting K\"unneth component in this case is
$H^1(X_1,\mQ) \otimes H^3(X_2,\mQ)$. As we have shown, $\Hg(X) = \{(u_1,u_2)
\in \UU_k \times \Hg(Y) \mid u_1^2 \cdot \det_k(u_2) = 1\}$. In particular we
have an element $(-1,1) \in \Hg(X)$ which acts on $H^1(X_1,\mQ) \otimes
H^3(X_2,\mQ)$ as $-1$. This shows there are no Hodge classes in $H^1(X_1,\mQ)
\otimes H^3(X_2,\mQ)$ and that $\cB^\gdot(X)$ is generated by divisor classes.
\whiteB\whiteB

\newcount\labno
\def\labelno{\number\labno}
\def\newlab{\advance\labno by 1}
\setbox0=\hbox{{\bf [99]}}
\def\mylabel{\newlab\item{[\labelno]}}
\def\refskip{\smallskip}
\def\author#1{{\rm #1}}
\def\title#1{{\sl #1}}

\frenchspacing
\eightpoint
\noindent
\centerline{\tenbf References.}
\whiteB

\noindent
\mylabel 
\author{M.~V. Borovoi},
\title{The Hodge group and the algebra of endomorphisms of an abelian variety},
In: A.~L. Onishchik (ed.), Problems in group theory and homological algebra,
Yaroslav. Gos. Univ., Yaroslavl (1981), 124--126. (MR 84m:14047).
\refskip

\noindent
\mylabel 
\author{N. Bourbaki},
\title{Groupes et Alg\`ebres de Lie},
Chap. 1, Hermann, Paris (1960); Chapters 7 et 8, Hermann, Paris (1975).
\refskip

\noindent
\mylabel 
\author{P. Deligne},
\title{La conjecture de Weil pour les surfaces K3},
Inventiones Math. {\bf 15} (1972), 206--226.
\refskip

\noindent
\mylabel 
\author{P. Deligne},
\title{Vari\'et\'es de Shimura: interpr\'etation modulaire, et techniques de
construction de mod\`eles canoniques},
In: A. Borel and W. Casselman (eds.), Automorphic forms, representations, and
$L$-functions, Proc. Symp. in pure Math. {\bf 33}(2), AMS, Providence (1979),
247--289.
\refskip

\noindent
\mylabel 
\author{F. Hazama},
\title{Hodge cycles on abelian varieties of CM-type},
Res. Act. Fac. Sci. Engrg. Tokyo Denki Univ. {\bf 5} (1983), 31--33. (MR
86g:14024).
\refskip

\noindent
\mylabel 
\author{F. Hazama},
\title{Algebraic cycles on certain abelian varieties and powers of special
surfaces},
J. Fac. Sci. Univ. Tokyo, Sect. \Romno 1a, {\bf 31} (1984), 487--520.
\refskip

\noindent
\mylabel 
\author{F. Hazama},
\title{Algebraic cycles on nonsimple abelian varieties},
Duke Math. J. {\bf 58} (1989), 31--37.
\refskip

\noindent
\mylabel 
\author{H. Imai},
\title{On the Hodge groups of some abelian varieties},
Kodai Math. Sem. Rep. {\bf 27} (1976), 367--372.
\refskip

\noindent
\mylabel 
\author{N. Jacobson},
\title{Lie Algebras},
Interscience Tracts {\bf 10}, John Wiley and Sons, New York (1962). Republished
by Dover Publications, Inc., New York (1979).
\refskip

\noindent
\mylabel 
\author{H.~W. Lenstra and Yu.~G. Zarhin},
\title{The Tate conjecture for almost ordinary abelian varieties over finite
fields},
In: F. Gouv\^ea and N. Yui (eds.), Advances in Number Theory: Proceedings of
the third conference of the CNTA, 1991, Clarendon Press, Oxford (1993),
179--194.
\refskip

\noindent
\mylabel 
\author{B.~J.~J. Moonen and Yu.~G. Zarhin},
\title{Hodge classes and Tate classes on simple abelian fourfolds},
Duke Math. J. {\bf 77} (1995), 553--581.
\refskip

\noindent
\mylabel 
\author{B.~J.~J. Moonen and Yu.~G. Zarhin},
\title{Weil classes on abelian varieties},
J. reine angew. Math. {\bf 496} (1998), 83--92.
\refskip

\noindent
\mylabel 
\author{D. Mumford},
\title{A Note of Shimura's paper ``Discontinuous groups and abelian
varieties''},
Math. Ann. {\bf 181} (1969), 345--351.
\refskip

\noindent
\mylabel 
\author{D. Mumford},
\title{Abelian varieties},
2nd Edition, Oxford Univ. Press, Oxford (1974).
\refskip

\noindent
\mylabel 
\author{V.~K. Murty},
\title{Exceptional Hodge classes on certain abelian varieties},
Math. Ann. {\bf268} (1984), 197--206.
\refskip

\noindent
\mylabel 
\author{F. Oort},
\title{Endomorphism algebras of abelian varieties},
In: H. Hijikata et al. (eds.), Algebraic geometry and commutative algebra in
honor of Masayoshi Nagata, Part \Romno 2, Kinokuniya Company, Ltd., Tokyo
(1988), 469--502.
\refskip

\noindent
\mylabel 
\author{R. Pink},
\title{$\ell$-adic algebraic monodromy groups, cocharacters, and the
Mumford-Tate conjecture},
J. reine angew. Math. {\bf 495} (1998), 187--237.
\refskip

\noindent
\mylabel 
\author{K.~Ribet},
\title{Hodge classes on certain types of abelian varieties},
Amer. J. Math. {\bf 105} (1983), 523--538.
\refskip

\noindent
\mylabel 
\author{C. Schoen},
\title{Hodge classes on self-products of a variety with an automorphism},
Compositio Math. {\bf 65} (1988), 3--32;
\title{Addendum}, Compositio Math. {\bf 114} (1998), 329--336.
\refskip

\noindent
\mylabel 
\author{J-P. Serre},
\title{Repr\'esentations lin\'eaires des groupes finis (troisi\`eme \'ed.)},
Hermann, Paris (1978).
\refskip

\noindent
\mylabel 
\author{J-P. Serre},
\title{Groupes alg\'ebriques associ\'es aux modules de Hodge-Tate},
In: Journ\'ees de g\'eometrie alg\'ebrique de Rennes, Ast\'erisque {\bf 65}
(1979), 155--188.
\refskip

\noindent
\mylabel 
\author{G. Shimura},
\title{On analytic families of polarized abelian varieties and automorphic
functions},
Ann. of Math. {\bf 78} (1963), 149--192.
\refskip

\noindent
\mylabel 
\author{S.~G. Tankeev},
\title{On Algebraic cycles on abelian varieties. \Romno 2},
Math. USSR Izv. {\bf 14} (1980), 383--394.
\refskip

\noindent
\mylabel 
\author{S.~G. Tankeev},
\title{Cycles on simple abelian varieties of prime dimension},
Math. USSR Izv. {\bf 20} (1983), 157--171.
\refskip

\noindent
\mylabel 
\author{B. van Geemen},
\title{Theta functions and cycles on some abelian fourfolds},
Math. Z. {\bf 221} (1996), 617--631.
\refskip

\noindent
\mylabel 
\author{A. Weil},
\title{Abelian varieties and the Hodge ring},
Collected papers, Vol. \Romno 3, [1977c], 421--429.
\refskip

\noindent
\mylabel 
\author{Yu.~G. Zarhin},
\title{Weights of simple Lie algebras in the cohomology of algebraic
varieties}, Math. USSR Izv. {\bf 24} (1985), 245--282.
\refskip

\noindent
\mylabel 
\author{Yu.~G. Zarhin},
\title{Linear irreducible Lie algebras and Hodge structures},
In: S. Bloch, I. Dolgachev, W. Fulton (eds.), Proceedings of the USA-USSR
Symposium on Algebraic
Geometry, Chicago 1989, Lecture Notes in Math. {\bf 1479}, Springer-Verlag,
Berlin (1991), 281--297.
\refskip

\noindent
\mylabel 
\author{Yu.~G. Zarhin},
\title{The Tate conjecture for non-simple abelian varieties over finite
fields},
In: G. Frey, J. Ritter (eds.), Algebra and Number Theory, de Gruyter, Berlin
(1994), 267--296.
\refskip

\noindent
\mylabel 
\author{Yu.~G. Zarhin and B.~J.~J. Moonen},
\title{Weil classes and {R}osati involutions on complex abelian varieties},
In: S.~G. Hahn, H.~C. Myung and E.~I. Zelmanov (eds.), Recent progress in
algebra (Taejon/Seoul, 1997), Contemp. Math. {\bf 224}, AMS, Providence (1999),
229--236.
\refskip
\bigskip

\noindent
B.~J.~J. Moonen, Department of Mathematics, University of Utrecht, P.O. Box
80.010, NL-3508 TA\ \ Utrecht, The Netherlands. E-mail: {\tt moonen@math.uu.nl}
\smallskip

\noindent
Yu.~G. Zarhin, Department of Mathematics, Pennsylvania State University,
University Park, PA 16802, USA. E-mail: {\tt zarhin@math.psu.edu}
\bye